\documentclass[11pt]{article}
\usepackage{bbm}
\usepackage{mathrsfs}
\usepackage{amsfonts}
\usepackage{amssymb}
\usepackage{amssymb,amsmath,amsthm, amsfonts}

\def\sqr#1#2{{\vcenter{\vbox{\hrule height.#2pt
              \hbox{\vrule width.#2pt height#1pt \kern#1pt \vrule width.#2pt}
          \hrule height.#2pt}}}}
\def\signed #1{{\unskip\nobreak\hfil\penalty50
          \hskip2em\hbox{}\nobreak\hfil#1
          \parfillskip=0pt \finalhyphendemerits=0 \par}}
\def\endpf{\signed {$\sqr69$}}
\def\sqr#1#2{{\vcenter{\vbox{\hrule height.#2pt
              \hbox{\vrule width.#2pt height#1pt \kern#1pt \vrule width.#2pt}
              \hrule height.#2pt}}}}
\def\signed #1{{\unskip\nobreak\hfil\penalty50
              \hskip2em\hbox{}\nobreak\hfil#1
              \parfillskip=0pt \finalhyphendemerits=0 \par}}
\def\endpf{\signed {$\sqr69$}}
\def\3n{\negthinspace \negthinspace \negthinspace }
\def\2n{\negthinspace \negthinspace }
\def\1n{\negthinspace }

\def\={\buildrel \triangle \over =}

\def\O{\Omega}

\def\q{\quad}

\def\max{\mathop{\rm max}}
\def\min{\mathop{\rm min}}

\def\sup{\mathop{\rm sup}}

\def\|{\Big |}
\def\({\Big (}
\def\){\Big )}
\def\[{\Big[}
\def\]{\Big]}
\def\be{\begin{equation}}
\def\bel{\begin{equation}\label}
\def\ee{\end{equation}}
\def\bt{\begin{theorem}}
\def\bcd{\begin{condition}}
\def\ecd{\end{condition}}
\def\et{\end{theorem}}
\def\bc{\begin{corollary}}
\def\ec{\end{corollary}}
\def\bde{\begin{definition}}
\def\ede{\end{definition}}
\def\bl{\begin{lemma}}
\def\el{\end{lemma}}
\def\bp{\begin{proposition}}
\def\ep{\end{proposition}}
\def\br{\begin{remark}}
\def\er{\end{remark}}
\def\ba{\begin{array}}
\def\ea{\end{array}}
\def\ed{\end{document}}

\def\square#1{\vbox{\hrule\hbox{\vrule height#1%
     \kern#1\vrule}\hrule}}
\def\rectangle#1#2{\vbox{\hrule\hbox{\vrule height#1%
     \kern#2\vrule}\hrule}}

\font\tenbb=msbm10 \font\sevenbb=msbm7 \font\fivebb=msbm5

\newfam\bbfam
\scriptscriptfont\bbfam=\fivebb \textfont\bbfam=\tenbb
\scriptfont\bbfam=\sevenbb

\newtheorem{lemma}{Lemma}[section]
\newtheorem{remark}{Remark}[section]
\newtheorem{example}{Example}[section]
\newtheorem{theorem}{Theorem}[section]
\newtheorem{corollary}{Corollary}[section]

\newtheorem{definition}{Definition}[section]
\newtheorem{proposition}{Proposition}[section]
\newtheorem{condition}{Condition}[section]

\makeatletter
   
   \@addtoreset{equation}{section}
\makeatother

\begin{document}

\title{Reflected Mean-Field Backward Stochastic Differential Equations. Approximation and Associated Nonlinear PDEs}

\author{Juan Li \bf\footnote{The author has
been supported by the NSF of P.R.China (Nos. 11071144, 11171187, 11222110), Shandong Province (Nos. BS2011SF010, JQ201202), 111 Project (No. B12023).}\\
{\small School of Mathematics and Statistics, Shandong University, Weihai, Weihai 264209, P. R. China.}\\
{\small {\it E-mail: juanli@sdu.edu.cn.}}}

\date{September 18, 2012}

\maketitle \noindent{\bf Abstract}

\hskip4mm Mathematical mean-field approaches have been used in many fields, not only in Physics and Chemistry, but also recently in Finance, Economics, and Game
Theory. In this paper we will study a new special mean-field problem in a purely probabilistic method, to characterize its limit which is the solution of mean-field backward stochastic differential equations (BSDEs) with reflections. On the other hand, we will prove that this type of reflected mean-field BSDEs can also be obtained as the limit equation of the mean-field BSDEs by penalization method. Finally, we give the probabilistic interpretation of the nonlinear and nonlocal partial differential equations with the obstacles by the solutions of reflected mean-field BSDEs.


\newpage
\section{\large{Introduction}}

\hskip0.5cm

Mathematical mean-field approaches have been used in many fields. To work on a stochastic limit approach to a mean-field problem is
inspired at the one hand by classical mean-field approaches in
Statistical Mechanics and Physics, by similar
methods in Quantum Mechanics and Quantum Chemistry, but also by a recent series of papers by Lasry and Lions
(see~\cite{LL} and the references inside cited) who studied mean-field games. And also it has been
strongly inspired by the McKean-Vlasov partial differential
equations (PDEs) which have found a great interest in the last years
and have been studied with the help of stochastic methods by many authors. On the other hand,
in the last years models of large stochastic particle systems with
mean field interaction have been studied by many authors; they have
described them by characterizing their asymptotic behavior when the
size of the system becomes very large, and also have shown that probabilistic methods
allow to study the solution of linear McKean-Vlasov PDE. The reader is referred, for
example, to the works by Borkar and Kumar~\cite{BK}, Bossy~\cite{B}, Bossy and Talay~\cite{BT},
Chan~\cite{C}, Kotelenez~\cite{K}, Mckean~\cite{MC}, M\'{e}l\'{e}ard~\cite{M},
Overbeck~\cite{O}, Pra and Hollander~\cite{PH},
Sznitman~\cite{S1},~\cite{S2}, Talay and Vaillant~\cite{TV}, and all
the references therein. More details may refer to Buckdahn, Djehiche, Li and
Peng~\cite{BLP1} and the references inside cited.

Buckdahn, Djehiche, Li and Peng~\cite{BLP1} studied a special mean-field
problem in a purely stochastic approach. They considered a stochastic
differential equation that describes the dynamics of a particle
$X^{(N)}$ influenced by the dynamics of $N$ other particles, which
are supposed to be independent identically distributed and of the
same law as $X^{(N)}$. This equation (of rank $N$) is then
associated with a backward stochastic differential equation (BSDE).
After having proven the existence and the uniqueness of a solution
$(X^{(N)},Y^{(N)},Z^{(N)})$ for this couple of equations the authors of~\cite{BLP1} investigated its limit behavior. With a new approach which uses the
tightness of the laws of the above sequence of triplets in a
suitable space, and combines it with BSDE methods and the Law of
Large Numbers, it was shown that $(X^{(N)},Y^{(N)},Z^{(N)})$
converges in $L^2$ to the unique solution of a limit equation formed
by a McKean-Vlasov stochastic differential equation and a Mean-Field
backward stochastic differential equation. Furthermore, Buckdahn, Li and Peng~\cite{BLP} proved the existence and the uniqueness of the solution of mean-field BSDEs under the classical assumptions, the comparison theorem of mean-field BSDEs and gave a stochastic interpretation to McKean-Vlasov partial differential equations (PDEs) with the help of the solutions of mean-field BSDEs. Since then we want to work on another new special mean-field problem to get a new limit equation which is like reflected BSDE in some sense. On the other hand, since the works~\cite{BLP} and~\cite{BLP1} on the mean-field BSDEs, there are many works on its generalizations, e.g., Wang~\cite{W} studied backward doubly SDEs of mean-field type and its applications; Shi, Wang and Yong~\cite{SWY} studied backward stochastic Volterra integral equations of mean-field type; Li and Luo~\cite{LL2} studied reflected BSDEs of mean-field type, they proved the existence and the uniqueness for reflected mean-field BSDEs; and also its applications, e.g.,  Andersson, Djehiche~\cite{AD}, Bensoussan, Sung, Yam and Yung~\cite{BSYY}, Buckdahn, Djehiche and Li~\cite{BBL}, Li~\cite{L}, Yong~\cite{Y}. Reflected BSDEs were introduced by El Karoui, Kapoudjian, Pardoux, Peng and Quenez~\cite{EKPPQ} in 1997. Later the theory of RBSDEs develops very quickly, because of its many applications, for example, in partial differential equations, finance and so on. More details may refer to Buckdahn and Li~\cite{BL} and the references inside cited.

In this paper we will study another new special mean-field problem, and get its limit which is a new type of reflected BSDEs, we call it reflected mean-field BSDEs. Our objective here is to characterize such an equation, at one hand, as the limit of classical BSDEs with
reflection and, on the other hand, as the limit of mean-field BSDEs
with a penalization approach. The approximating reflected BSDEs (N) are
discussed, and with an example it is in particular shown that these
reflected BSDEs (N)  don't obey the comparison principle. Furthermore,
under an additional monotonicity assumption of the driving
coefficient, the description of reflected mean-field BSDEs as
monotonic limit of mean-field BSDEs without reflection is used to give
through them a stochastic interpretation of associated non-local PDEs
with obstacles. We show that the solution of the reflected mean-field
BSDE is the unique viscosity solution of the associated non-local
PDE with obstacles.

More precisely, we consider the following mean-field BSDE with reflections:
\be\label{1}
  \begin{array}{ll}
&{\rm (i)}\ Y \in S_{\mathbf{F}}^2([0,T]), \, Z \in
L_{\mathbf{F}}^2([0,T];{\mathbb{R}}^d)\ \mbox{and}\  K\in
A_{\mathbf{F}}^{2,c}([0,T]);\\
& {\rm (ii)}\ Y_t = \displaystyle E\left[\Phi(x,X_T)\right]_{x=X_T} +
\int_t^TE\left[g(s,\mathbf{u},\Lambda_s)\right]_{\vert
\mathbf{u}=\Lambda_s}ds + K_{T}-K_{t} -\int^T_tZ_sdW_s;\\
& {\rm (iii)}\ \displaystyle Y_t \geq h(t,X_t),\ \mbox{a.s., for all}\  t\in
[0,T];\\
& {\rm (iv)}\ \displaystyle\int_0^T(Y_t - h(t,X_t))dK_{t}=0,\\
\end{array}
\ee
\noindent where we have used the notation $\Lambda=(X,Y,Z);$ $T>0$\ is a given finite time horizon; $W=(W_t)_{t\geq 0}$\ is a d-dimensional Brownian motion; $X=(X_t)_{t\in [0, T]}$\ is a driving n-dimensional adapted stochastic process.

Such type of mean-field BSDEs without reflections have been studied by Buckdahn, Li and Peng~\cite{BLP}, they proved that such a mean-field BSDE gave a stochastic interpretation to the related nonlocal PDEs. In this paper we first prove that, under our standard assumptions the mean-field BSDE (\ref{1}) with reflections will be the limit equation of the following reflected BSDE (N):
\be
  \begin{array}{ll}
& {\rm (i)}\ Y^N \in S_{\mathbf{F}}^2([0,T]), \, Z^N \in
L_{\mathbf{F}}^2([0,T];{\mathbb{R}}^d)\ \mbox{and}\ K^N\in
A_{\mathbf{F}}^{2,c}([0,T]);\\
& {\rm (ii)}\ Y^N_t = \displaystyle \xi^N +
\int_t^Tf^N(s,\Theta_N(Y_s,Z_s))ds + K^N_{T}-K^N_{t}
-\int^T_tZ^N_sdW_s,\ t\in [0,T];\\
& {\rm (iii)}\ Y^N_t \geq L^N_t,\ \mbox{a.s., for any}\  t\in [0,T];\\
& {\rm (iv)}\ \displaystyle\int_0^T(Y^N_t - L^N_t)dK^N_{t}=0,
\end{array}
\ee
where for $N\ge 1$ and $\omega\in\Omega$,
$$
  \begin{array}{cl}
& \displaystyle \xi^N(\omega):=\frac{1}{N}\sum_{k=1}^N\Phi(\Theta^k(\omega),X^N_T(\omega),X^N_T(\Theta^k(\omega))),\\
& \displaystyle f^N(\omega,t,\textbf{y},\textbf{z}):= \frac{1}{N}\sum_{k=1}^N g(\Theta^k(\omega),t,X^N_t(\omega),(y_0,z_0),X^N_t(\Theta^k(\omega)),(y_k,z_k)),\\
&\mbox{\quad for }t\in[0,T],\, \textbf{y}=(y_0,\cdots,y_N)\in {\mathbb{R}}^{N+1},\textbf{z}=(z_0,\cdots,z_N)\in {\mathbb{R}}^{(N+1)\times d},\\
& \\
& \displaystyle L^N_t(\omega):=h(\omega,t,X^N_t(\omega)),\,
t\in[0,T].
\end{array}
$$
(More details refer to Theorem 4.1). Example 3.1 shows that such reflected BSDE (N) usually doesn't have the comparison theorem.

Furthermore, more generally, for the obstacle process $L=(L_s)_{0\leq s\leq T}\in S_{\mathbf{F}}^2([0,T])$\ and the terminal condition $\xi\in L^2(\Omega, {\cal F}_T, \mathbb{R})$ such that $\xi\geq L_T$, P-a.s., we consider the following reflected mean-field BSDE:
\be\label{3}
\begin{array}{rcl}
& & {\rm (i)} Y \in S_{\mathbf{F}}^2([0,T]), \, Z \in
L_{\mathbf{F}}^2([0,T];{\mathbb{R}}^d)\ \mbox{and}\ K\in
A_{\mathbf{F}}^{2,c}([0,T]);\\
& & {\rm (ii)} Y_t = \displaystyle \xi + \int_t^TE\left[g(s,y,z,Y_s)\right]_{\vert
y=Y_s, z=Z_s}ds + K_{T}-K_{t} -\int^T_t Z_sdW_s;\\
& & {\rm (iii)} \displaystyle Y_t \geq L_t,\ \mbox{a.s.,\ for all}\ t\in
[0,T]; \\
& & {\rm (iv)} \displaystyle\int_0^T(Y_t - L_t)dK_{t}=0,
\end{array}
\ee
Under our assumptions the reflected mean-field BSDE (\ref{3}) can be also got as the limit equation of the following penalized mean-field BSDEs:
\be\label{4} \widetilde{{Y}}_t^n=\xi+\int_t^T E[g(s, y, z, \widetilde{Y}^n_s)]|_{y=\widetilde{Y}^n_s, z=\widetilde{Z}^n_s}ds+n\int_t^T(\widetilde{{Y}}_s^n-L_s)^{-}ds-\int_t^T\widetilde{{Z}}_s^ndW_s.\ee
(More details refer to Theorem 5.1).

Finally, this allows us to give a probabilistic representation for the solution of the following non-local PDEs with obstacles:
\be\label{6.5}\left
\{\begin{array}{ll}
 &\!\!\!\!\! min\{u(t,x)-h(t,x),-\frac{\partial }{\partial t} u(t,x)-Au(t,x)\\
&\hskip2cm-E[f(t,x, X_t^{0, x_0},  u(t, X_t^{0, x_0}), u(t,x),
Du(t,x).E[\sigma(t, x, X_t^{0, x_0})])]\}=0,\\
 & \hskip 4.3cm (t,x)\in [0,T)\times {\mathbb{R}}^n ,  \\
 &\!\!\!\!\!  u(T,x) =E[\Phi (x, X_T^{0, x_0})], \hskip0.5cm   x \in
 {\mathbb{R}}^n,
 \end{array}\right.
\ee with
$$ Au(t,x):=\frac{1}{2}tr(E[\sigma(t,x, X_t^{0, x_0})]
 E[\sigma(t, x, X_t^{0, x_0})]^{T}D^2u(t,x))+Du(t,x).E[b(t,
 x, X_t^{0, x_0})].$$ Here the functions $b, \sigma, f\ \mbox{and}\ \Phi$\ are
supposed to satisfy (H6.1), and (H6.2), respectively, and
$X^{0, x_0}$\ is the solution of the SDE (\ref{6.1}). More details refer to Theorem 6.1.

Our paper is organized as follows: Section 2 recalls briefly some
elements of the theory of backward SDEs and mean-field BSDEs which will be needed in what
follows. In Section 3 we introduce the reflected BSDEs of rank $N$,
define the framework in which it is investigated and prove the
existence and the uniqueness, and give an example to explain that this type of reflected BSDEs of rank $N$ usually doesn't have the comparison theorem anymore. In this section we also give an important inequality about RBSDE which is very useful-Lemma 3.5. In Section 4 we prove the convergence of the solution of the reflected BSDE of rank (N) to that of reflected  mean-field BSDE (Theorem 4.1). In Section 5 we prove that the reflected mean-field BSDEs can also be obtained as the limit equations of the reflected BSDEs with the help of the penalization method (Theorem 5.1).  In Section 6 we prove that the solution of the reflected mean-field BSDEs is the unique viscosity solution of the associated nonlinear and nonlocal partial differential equation with the obstacles (Theorem 6.1). We also prove that the value functions $u_n(t,x)$ which are defined by the penalized mean-field BSDEs are Lipschitz with respect to $x$, uniformly in $t$, and $n\in \mathbb{N}$ (Proposition 6.1).

\section{ {\large Preliminaries}}

\subsection{A Recall on BSDEs}

In this section we will introduce some basic notations and
results about BSDEs, which will be needed in the following sections. First we will extend slight the classical Wiener space $(\Omega,{\cal F},P):$

\smallskip

$\bullet$ For an arbitrarily given time horizon $T>0$ and a
countable index set $I$ (which will be clarified later), $\Omega$ is
the set of all families $(\omega^i)_{i\in I}$, where $\omega^i:[0,T]\rightarrow{\mathbb{R}}^d$ is continuous with initial value $0$ (i.e., $\Omega=C_0([0,T];{\mathbb{R}}^d)^I$); we endow it with the product topology produced by the uniform convergence on its components
$C_0([0,T];{\mathbb{R}}^d);$

$\bullet$ Let ${\cal B}(\Omega)$ denote the Borel $\sigma$-field over
$\Omega$ and $B=(W^i)_{i\in I}$ be the coordinate process over
$\Omega:$ $W^i_t(\omega)=\omega^i_t,\,
t\in[0,T],\omega\in\Omega,i\in I;$

$\bullet$ Let $P$ be the Wiener measure over $(\Omega,{\cal
B}(\Omega)),$ i.e., the coordinates $W^i,\, i\in I,$ are a family of independent
$d$-dimensional Brownian motions with respect to $P$. In the end,

$\bullet$ Let ${\cal F}$ be the $\sigma$-field ${\cal B}(\Omega)$
completed by the Wiener measure $P$.

\smallskip

\noindent Define $W:=W^0.$ The Probability space $(\Omega,{\cal
F},P)$ is endowed with the filtration $\mathbf{F}=({\cal F}_t)_{t\in[0,T]}$, where $\mathbf{F}$ is generated by the Brownian motion $W$, enlarged by the
$\sigma$-field ${\cal G}=\sigma\{W^i_t,\, t\in[0,T],i\in
I\setminus\{0\}\}$ and completed by the collection ${\cal N}_P$ of
all $P$-null sets, that is
$${\cal F}_t={\cal F}^W_t\vee{\cal G},\, \, t\in[0,T],$$
where $\mathbf{F}^W=({\cal F}^W_t=\sigma\{W_r,\, r\le
t\}\vee{\cal N}_P)_{t\in[0,T]}.$ Notice that the Brownian motion $W$ still has the martingale representation
property with respect to the filtration $\mathbf{F}$ .

We also introduce the following spaces which will be used later:
$$S_{\mathbf{F}}^2([0,T])=\{(Y_t)_{t\in [0,T]} \mbox{ continuous adapted
process:}\ E[\sup_{t\in[0,T]}\vert Y_t\vert^2]
<+\infty\};$$
$$\begin{array}{rcl}
L_{\mathbf{F}}^2([0,T];{\mathbb{R}}^d) & = & \{(Z_t)_{t\in [0,T]}
\, {\mathbb{R}}^d\mbox{-valued progressively measurable process: }\\
& &\hskip 2cm E\left[\int_{0}^T\vert Z_t\vert^2dt\right]
<+\infty\}.\end{array}$$

\medskip

\noindent (Notice that $\vert z\vert$ denotes the Euclidean norm of
$z\in {\mathbb{R}}^d$).  Now given a measurable function
$g:\Omega\times[0,T]\times {\mathbb{R}}\times
{\mathbb{R}}^d\rightarrow {\mathbb{R}} $\ which satisfies that
$(g(t, y,z))_{t\in [0, T]}$ is ${\mathbf{F}}$-progressively
measurable for all $(y,z)$ in ${\mathbb{R}}\times {\mathbb{R}}^d$.
We give the following standard assumptions:
\vskip0.2cm

(A1) There is some constant $C\ge 0$  such that, P-a.s., for all $t\in
[0, T],\
y_1, y_2\in {\mathbb{R}},\ z_1, z_2\in {\mathbb{R}}^d,\\
\mbox{ } \hskip1.9cm   \vert g(t, y_1, z_1) - g(t, y_2,
z_2)\vert\leq C(|y_{1}-y_{2}| + |z_{1}-z_{2}|).$
 \vskip0.2cm

(A2) $g(\cdot,0,0)\in L_{\mathbf{F}}^2([0,T];{\mathbb{R}})$.
\vskip0.2cm

The following results on BSDEs are classical; for the proof refer to, e.g., Pardoux and Peng~\cite{PaPe}, or El Karoui, Peng and
Quenez~\cite{ElPeQu}.

\medskip

\bl Suppose the generator $g$ satisfies (A1) and (A2). Then, for any random variable $\xi\in
L^2(\Omega,{\cal F}_T,P),$ the BSDE associated with the data $(g, \xi)$
\be Y_t = \xi + \int_t^Tg(s,Y_s,Z_s)ds - \int^T_tZ_s\, dB_s,\, \, 0\le t\le T,\ee
has a unique solution $(Y, Z)\in {\cal{S}_{\mathbf{F}}}^2([0, T])\times L_{\mathbf{F}}^2([0,T];{\mathbb{R}}^d).$\el

\medskip

Now we give the standard estimate for BSDEs.

\medskip

\bl Suppose that $g_k$ satisfies (A1) and (A2) and $\xi_k\in L^{2}(\Omega,
{\cal{F}}_{T}, P),\, k=1,2.$ Let $(Y^k,Z^k)$ denote the unique
solution of the BSDE with the data $(g_k,\xi_k),\, k=1,2,$\ respectively. For every
$\delta>0$, there exists some $\gamma>0, C>0$ only depending on $\delta$ and the Lipschitz
constants of $g_k,k=1,2,$ such that,
$$\begin{array}{cl}
  & \displaystyle E[\int_0^T e^{\gamma t}(\vert\overline{Y}_t\vert^2+\vert\overline{Z}_t \vert^2)dt]  \le\displaystyle CE[e^{\gamma T}\vert\overline{\xi}\vert^2]+ \delta E[\int_0^T e^{\gamma t}\vert\overline{g}(t,Y^1_t,Z^1_t)\vert^2dt],
  \end{array}$$
where
$$(\overline{Y},\overline{Z})=(Y^1-Y^2,Z^1-Z^2),\, \,
\overline{g}=g_1-g_2,\, \, \overline{\xi}=\xi_1-\xi_2.$$\el

\medskip

Now we introduce one of the important results for BSDEs-the
comparison theorem (see Proposition 2.4 in Peng~\cite{P} or Theorem 2.2 in El Karoui, Peng and
Quenez~\cite{ElPeQu}).

\medskip

\bl (Comparison Theorem) Suppose two
coefficients $g_1$ and $g_2$ satisfy (A1) and (A2) and two
terminal values $ \xi_1,\ \xi_2 \in L^{2}(\Omega, {\cal{F}}_{T},
P)$. $(Y^1,Z^1)$\ and $(Y^2,Z^2)$\ are the solutions of the
BSDE with the data $(\xi_1,g_1 )$\ and $(\xi_2,g_2 )$, respectively.
Then we have:

{\rm (i) }(Monotonicity) If  $ \xi_1 \geq \xi_2$  and $ g_1 \geq
g_2, \ a.s.$, then $Y^1_t\geq Y^2_t$, for all $t\in [0, T]$, a.s.

{\rm (ii)}(Strict Monotonicity) If, in addition to {\rm (i)},  also
$P\{\xi_1 > \xi_2\}> 0$, then we have $P\{Y^1_t> Y^2_t\}>0,$ for all
$\ 0 \leq t \leq T,$\ and in particular, $Y^1_0> Y^2_0.$\el

\vskip0.5cm

\subsection{A Recall on Mean-field BSDEs}

\vskip0.5cm

In this section we recall some basic results on a new type of BSDE, the so
called Mean-Field BSDEs; more details refer to~\cite{BLP1} and~\cite{BLP}.

\noindent The driver of our mean-field BSDE is a function
$f=f(\omega, t, y, z, \tilde{y},\tilde{z}):  {\Omega}\times[0,T]\times
{\mathbb{R}} \times {\mathbb{R}}^{d}\times {\mathbb{R}} \times
{\mathbb{R}}^{d} \rightarrow {\mathbb{R}}$\ which is
${{\mathbf{F}}}$-progressively measurable, for all
$(y, z, \tilde{y}, \tilde{z})$, and satisfies the following assumptions:
 \vskip0.2cm

(A3) There exists a constant $C\ge 0$  such that, ${P}$-a.s.,
for all $t\in [0, T],\ y_{1}, y_{2}, \tilde{y}_{1}, \tilde{y}_{2}\in
{\mathbb{R}},\ z_{1}, z_{2},\tilde{ z}_{1}, \tilde{z}_{2}\in {\mathbb{R}}^d,\\
\mbox{ }\hskip1.5cm |f(t, y_{1}, z_{1}, \tilde{y}_{1}, \tilde{z}_{1}) - f(t, y_{2}, z_{2},
\tilde{y}_{2}, \tilde{z}_{2})|\leq C(|y_{1}-y_{2}|+ |z_{1}-z_{2}|+|\tilde{y}_{1}-\tilde{y}_{2}|+
|\tilde{z}_{1}-\tilde{z}_{2}|).$
 \vskip0.2cm

(A4) $f(\cdot,0,0,0,0)\in
{L}_{{{\mathbf{F}}}}^{2}(0,T;{\mathbb{R}})$. \vskip0.2cm

The following results refer to~\cite{BLP}.

\bl Under the assumptions (A3) and (A4), for any random variable
$\xi\in L^2(\O, {\cal{F}}_T,$ $P),$ the mean-field BSDE
 \be Y_t = \xi + \int_t^TE[f(s,y,z,Y_s,Z_s)]|_{y=Y_s, z=Z_s}ds - \int^T_tZ_s\,
dB_s,\q 0\le t\le T, \label{BSDE1} \ee
 has a unique adapted solution
$$(Y_t, Z_t)_{t\in [0, T]}\in {\cal{S}}_{{\mathbf{F}}}^2(0, T; {\mathbb{R}})\times
{L}_{{\mathbf{F}}}^{2}(0,T;{\mathbb{R}}^{d}). $$ \el

\bl (Comparison Theorem) Let $f_i=f_i({\omega}, t, y, z, \tilde{y},\tilde{z}),\ i=1, 2,$\ be two drivers satisfying the standard assumptions
(A3) and (A4). Moreover, we suppose

\noindent{\rm(i)} One of both coefficients is independent of $\tilde{z}$;

\noindent{\rm(ii)} One of both coefficients is nondecreasing in
$\tilde{y}$.

Let $\xi_1,\ \xi_2\in L^{2}(\Omega, {\cal{F}}_{T}, P)$\ and denote
by $(Y^1, Z^1)$\ and $(Y^2, Z^2)$\ the solution of the mean-field
BSDE (\ref{BSDE1}) with data $(\xi_1, f_1)$\ and $(\xi_2, f_2)$,
respectively. Then if  $\xi_1\leq \xi_2,\ \mbox{P-a.s.}$, and
$f_1\leq f_2,\  {P}\mbox{-a.s.},$\ it holds that also $Y_t^1\leq
Y_t^2,\ t\in [0, T],\ \mbox{P-a.s.}$\el

\br The conditions {\rm(i)} and {\rm(ii)} of Lemma 2.2 are, in
particular, satisfied, if they hold for the same driver $f_j$\ but
also if {\rm(i)} is satisfied by one driver and {\rm(ii)} by the
other one. \er

\section{Reflected BSDEs}

After the short recall on BSDEs let us now consider
reflected BSDEs (RBSDEs) and MFBSDEs with reflection. Let us first
introduce the framework in which we want to study the limit approach
to get reflected MFBSDEs. First we give the countable index set as follows:
$$I:=\{i\, \vert\, i\in\{1,2,3,\dots\}^ k,\, k\ge 1\}\cup \{0\}.$$
We define $i\oplus i'=(i_1,\dots,i_k,i'_1,\dots,i'_{k'})\in I$, for two elements $i=(i_1,\dots,i_k),$ $i'=(i'_1,\dots,i'_{k'})$ of
$I$ (with the convention that $i\oplus 0=i$); in particular, $\omega^{(k)\oplus i}=\omega^{(k, i_1,\cdots, i_k)}, k\geq 0$. Notice that $\left(\omega^{(k)\oplus i}\right)_{i\in I}\in \Omega$.

Now we introduce a family of shift operators
$\Theta^k:\Omega\rightarrow\Omega$, $k\ge 0,$ over $\Omega$. We define $\Theta^k(\omega)=(\omega^{(k)\oplus i})_{i\in I},\,
\omega\in\Omega,k\ge 0,$ and notice that $\Theta^k$ is an operator
mapping $\Omega$ into $\Omega$\ associating $(\omega^i)_{i\in I}\in \Omega$\ with $(\omega^{(k)\oplus i})_{i\in I}\in \Omega$. Notice that all these operators
$\Theta^k:\Omega\rightarrow\Omega$ make the Wiener measure $P$
invariant (i.e., $P\circ[\Theta^k]^{-1}=P$), which allows to regard
$\Theta^k$ as an operator defined over $L^0(\Omega,{\cal F},P):$
\ $\Theta^k(\xi)(\omega):=\xi(\Theta^k(\omega)),\,
\omega\in\Omega,$ for the random variables $\xi(\omega)=f(\omega^{i_1}_{t_1}, \cdots,\omega^{i_n}_{t_n}),$
$i_1,\cdots,i_n\in I,\, t_1,\cdots,t_n\in[0,T],\, f\in C(R^{d\times
n}),\, n\ge 1,$ then we can  extend this definition from this set of
continuous Wiener functionals to the space $L^0(\Omega,{\cal
F},P)$ with the help of the density of the set of smooth Wiener
functionals in $L^0(\Omega,{\cal F},P)$. Notice that, for all
$\xi\in L^0(\Omega,{\cal F},P),$ the random variables
$\Theta^k(\xi),\, k\ge 1,$ are independent and uniformly identically distributed
(i.i.d.), with the same law as $\xi$ and also independent of the Brownian
motion $W$.

In the end, for simplicity of notations we introduce the $(N+1)$-dimensional
shift operator $\Theta_N=(\Theta^0,\Theta^1,\cdots,\Theta^N)$, which
relates a random variable $\xi\in L^0(\Omega,{\cal F},P)$ with
the $(N+1)$-dimensional random vector
$\Theta_N(\xi)=(\xi,\Theta^1(\xi),\cdots,$ $\Theta^N(\xi))$ (remark
that $\Theta^0$ is the identical operator). If $\xi$ is
a random vector, $\Theta^k(\xi)$ and $\Theta_N(\xi)$ are introduced by
a componentwise application of the corresponding operators.

\smallskip

For introducing the notion of a RBSDE we shall introduce still the
following space of adapted increasing processes:
$$A_{\mathbf{F}}^{2,c}([0,T])=\{(K_t)_{t\in [0,T]}\in S_{\mathbf{F}}^2([0,T])\mbox{ non-decreasing process: }K_0=0\}.$$
An RBSDE with one barrier is associated with a terminal condition
$\xi \in L^{2}(\Omega,{\cal{F}}_{T},$ $P)$, a generator $g$
satisfying the assumptions (A1) and (A2), and an ``obstacle" process
$L=(L_t)_{t\in[0,T]}$. We assume that $L\in S^2_{\mathbf{F}}([0,T])$
and $L_T\leq\xi, \ P$-a.s. A solution of an RBSDE with one barrier
is a triplet $(Y, Z, K)$ of $\mathbf{F}$-progressively measurable
processes, taking its values in ${\mathbb{R}}\times
{\mathbb{R}}^d\times {\mathbb{R}}_+$ and satisfying the following
properties

\be
\begin{array}{lll}
& {\rm (i)}\ Y \in S_{\mathbf{F}}^2([0,T]), \, Z \in
L_{\mathbf{F}}^2([0,T];{\mathbb{R}}^d)\ \mbox{and}\ K\in A_{\mathbf{F}}
^{2,c}([0,T]);\\
& {\rm (ii)}\ Y_t = \displaystyle \xi + \int_t^Tg(s,Y_s,Z_s)ds +
K_{T}-K_{t} -\int^T_tZ_sdW_s,\quad t\in [0,T]; \\
& {\rm (iii)}\ Y_t \geq L_t,\ \mbox{a.s., for any}\ t\in [0,T];\\
& {\rm (iv)}\ \displaystyle\int_0^T(Y_t - L_t)dK_{t}=0.
\end{array}\ee

\smallskip

\noindent For shortness, a given triplet $(\xi, g, L)$ is said to
satisfy the Standard Assumptions (A) if the generator $g$ satisfies
(A1) and (A2), the terminal value $\xi$ belongs to $L^2(\Omega,{\cal
F}_T,P)$, and the obstacle process $L\in S_{\mathbf{F}}^2([0,T])$ is
such that $L_T\le\xi$, $P$-a.s.

We begin by recalling two lemmata which are by now well-known results of
the theory of reflected BSDEs and are borrowed from Theorem 5.2 and
Theorem 4.1, respectively, of the paper by El Karoui, Kapoudjian,
Pardoux, Peng, Quenez~\cite{EKPPQ}.

\bl Let $(\xi,g,L)$ be a triplet satisfying
the Standard Assumptions (A). Then the above RBSDE admits a unique
solution $(Y, Z, K).$\el
\bl (Comparison Theorem) We suppose that two
triplets $(\xi_1, g_1,L^1)$ and $(\xi_2, g_2, L^2)$\ satisfy the
Standard Assumptions (A) but we impose only for one of the both
coefficients $g_1$\ and $ g_2$ to fulfill the Lipschitz condition
(A1). Furthermore, we make the following assumptions:
\be
  \begin{array}{ll}
{\rm(i)}&\xi_1 \leq \xi_2,\ \ a.s.;\\
{\rm(ii)}&g_1(t,y,z) \leq g_2(t,y,z),\ a.s., \hbox{ \it for all }
(t,y,z)\in [0,T]\times
{\mathbb{R}}\times {\mathbb{R}}^d;\\
{\rm(iii)}& L_t^1 \leq L^2_t,\ \ a.s., \hbox{ \it for all } t\in [0,T]. \\
 \end{array}
 \ee
Let $(Y^1,Z^1, K^1)$ and $(Y^2, Z^2, K^2)$ be adapted solutions of the RBSDEs
with data $(\xi_1, g_1,L^1)$ and $(\xi_2, g_2, L^2),$ respectively.  Then,
$Y^1_{t} \leq Y^2_{t}, $ for all $t\in [0,T],$ $P$-a.s.\el

\bigskip

We also shall recall the following both standard estimates of BSDEs
with one reflecting barrier.

\medskip

\bl Let $(Y,Z,K)$ be the solution of the
above RBSDE with data $(\xi,g,L)$\ satisfying the Standard
Assumptions (A). Then there exists a constant $C$\ such that
\be
  \begin{array}{ll}
 & \displaystyle E [\sup_{t\leq s\leq T }|Y_s|^2+\int_t^T|Z_s|^2ds+|K_T-K_t|^2 |{{\cal {F}}_t} ]\\
\leq & \displaystyle CE[\xi^2+\left(\int_t^T g(s,0,0)ds\right)^2+\sup_{t\leq s\leq
T} L_s^2|{{\cal {F}}_t} ],\, P\mbox{-a.s., } t\in[0,T]. \\
\end{array}
 \ee
The constant $C$\ depends only on the Lipschitz constant of $g$.
\el

Lemma 3.3 is based on Propositions 3.5 in~\cite{EKPPQ} and its
generalization by Proposition 2.1 in Wu and Yu~\cite{WY}. The
following statement refers to Proposition 3.6 in~\cite{EKPPQ}
or Proposition 2.2 in~\cite{WY}.

\medskip

\bl Let $(\xi,g,L)$ and $(\xi',g',L')$ be
two triplets satisfying the above Standard Assumptions (A). We
suppose that $(Y,Z,K)$ and $(Y',Z',K')$ are the solutions of our
RBSDE with the data $(\xi,g,L)$ and $(\xi',g',L'),$\ respectively.
Then, for some constant $C$ which only depends on the Lipschitz
constant of the coefficient $g'$, and with the notations
$$
  \begin{array}{cc}
& \overline{\xi}=\xi-\xi',\ \ \  \overline{g}=g-g',\ \ \  \overline{L}=L-L',\\
& \overline{Y}=Y-Y',\ \ \  \overline{Z}=Z-Z',\ \ \  \overline{K}=K-K',\\
  \end{array}
$$
it holds, for all $t\in[0,T],$ $P$-a.s.,
\be
  \begin{array}{cl}
& \displaystyle E[\sup_{s\in[t,T]}\vert\overline{Y}_s\vert^2+\int_t^T\vert\overline {Z}_s\vert^2ds+\vert\overline{K}_T-\overline{K}_t\vert^2\vert{\cal F}_t]\\
& \displaystyle \le C E[\vert\overline{\xi}\vert^2+\left(\int_t^T\vert\overline{g}(s,Y_s,Z_s)\vert ds\right)^2\vert{\cal F}_t]+C \left(E[\sup_{s\in[t,T]}\vert\overline{L}_s\vert^2\vert{\cal F}_t]\right)^{1/2} \Psi_{t,T}^{1/2},\\
  \end{array}
\ee
where
$$
  \begin{array}{cl}
& \displaystyle\Psi_{t,T}  =  E[\vert\xi\vert^2+\left(\int_t^T\vert g(s,0,0)\vert ds\right)^2+\sup_{s\in[t,T]} \vert L_s\vert^2\\
& \qquad\, \, \displaystyle +\vert\xi'\vert^2+\left(\int_t^T\vert g'(s,0,0)\vert ds\right)^2+\sup_{s\in[t,T]}\vert L'_s \vert^2\vert{\cal F}_t].\\
\end{array}
$$
\el

However we will also need a slight version of the above standard
estimate for RBSDEs, which is of the same nature as that given by
Lemma 2.2 for BSDEs.

\medskip

\bl As in Lemma 3.4 we suppose that
$(\xi,g,L)$ and $(\xi',g',L')$ are two triplets satisfying the above
Standard Assumptions (A), and $(Y,Z,K)$ and $(Y',Z',K')$ are the
solutions of the associated RBSDEs, respectively. Then, for some $C>0$ and for all $\delta>0$
we can find $\gamma>0$ such that, with
the notations of Lemma 3.4,
\be
  \begin{array}{cl}
\displaystyle E\left[\int_0^Te^{\gamma t}(\vert\overline{Y}_t\vert^2+\vert
\overline{Z}_t\vert^2) dt\right]\le & \displaystyle 2E\left[e^{\gamma T}\vert
\overline{\xi}\vert^2\right]+\delta E\left[ \int_0^Te^{\gamma t} \vert
\overline{g}(t,Y_t,Z_t)\vert^2dt \right]\\
& +\displaystyle C\left(E\left[\sup_{t\in[0,T]}\vert e^{\gamma
t}\overline{L}_t\vert^2\right]\right)^{1/2} \Psi_{0,T}^{1/2}
\end{array}
\ee
(Recall the definition of $\Psi_{t,T}$ given in Lemma 3.4.). The
constant $C$ only depends on the bound and the Lipschitz constant of
$g$ and $g'$, while $\gamma$ only depends on $\delta$ and on the Lipschitz
constant of $g$ and $g'$.\el
\br The above estimate for $L=L'$ was
established in the proof of Theorem 5.2 in~\cite{EKPPQ}. Since the
above lemma plays a crucial role in our approach we give its proof
for the reader's convenience.\er

\noindent\textbf{Proof} (of Lemma 3.5). Let $\delta>0$ be
sufficiently small and $\gamma>0$. Then, by applying It\^{o}'s
formula to the process $(e^{\gamma
t}\vert\overline{Y}_t\vert^2)_{t\in[0,T]}$ and by taking into
account that $(\overline{Y}_t-\overline{L}_t)d\overline{K}_t\le 0,$
$t\in[0,T],$ we get
$$
  \begin{array}{cl}
& \displaystyle E\left[\int_0^Te^{\gamma t}(\gamma\vert\overline{Y}_t\vert^2+\vert \overline{Z}_t\vert^2) dt\right]\\
& \le E\left[e^{\gamma T}\vert\overline{\xi}\vert^2\right]+\displaystyle 2E\left[ \int_0^Te^{\gamma t}\overline{Y}_t(g'(t,Y_t,Z_t)- g'(t,Y'_t,Z'_t))dt\right]\\
&\mbox{ }\hskip0.5cm + \displaystyle 2E\left[ \int_0^Te^{\gamma
t}\overline{Y}_t\overline{g}(t,Y_t,Z_t)dt
\right]+2E\left[\int_0^Te^{\gamma
t}\overline{L}_td\overline{K}_t\right].
\end{array}
$$
Thus, since $g'(\omega,t,.,.)$ is Lipschitz, with a Lipschitz
constant $L$ which does not depend on $(\omega,t)$, we can conclude
from the latter relation that, for some constant $C_{L,\delta}$ only
depending on $\delta$ and $L$,
$$
  \begin{array}{cl}
& \displaystyle E\left[\int_0^Te^{\gamma t}(\gamma\vert\overline{Y}_t\vert^2+\vert \overline{Z}_t\vert^2) dt\right]\\
& \le E\left[e^{\gamma T}\vert\overline{\xi}\vert^2\right]+ \displaystyle E\left[ \int_0^T e^{\gamma t}(C_{L,\delta}\vert\overline{Y}_t\vert^2+\frac{1}{2}\vert \overline{Z}_t\vert^2) dt\right]\\
&\mbox{ }\hskip0.5cm + \displaystyle \frac{1}{2}\delta E\left[
\int_0^Te^{\gamma t}\vert\overline{g}(t,Y_t,Z_t)\vert^2dt
\right]+2E\left[\sup_{t\in[0,T]}\vert e^{\gamma
t}\overline{L}_t\vert (K_T+K'_T)\right].
\end{array}
$$
Consequently, for $\gamma:=C_{L,\delta}+\frac{1}{2}$,
\be
  \begin{array}{cl}
& \displaystyle E\left[\int_0^Te^{\gamma t}(\vert\overline{Y}_t\vert^2+\vert \overline{Z}_t\vert^2) dt\right]\le 2E\left[e^{\gamma T}\vert\overline{\xi}\vert^2\right]+\delta E\left[ \int_0^Te^{\gamma t}\vert\overline{g}(t,Y_t,Z_t)\vert^2dt \right]\\
&\mbox{ }\hskip 0.5cm +\displaystyle
4\left(E\left[\sup_{t\in[0,T]}\vert e^{\gamma
t}\overline{L}_t\vert^2\right]\right)^{1/2}
\left(E\left[(K_T+K'_T)^2\right]\right)^{1/2}.
\end{array}
\ee
Then the result announced in the lemma follows from Lemma 3.3.\endpf

For an arbitrarily given natural number $N\ge 0$ we consider a
measurable function $f:\Omega\times [0,T]\times
{\mathbb{R}}^{N+1}\times {\mathbb{R}}^{(N+1)\times d}\rightarrow
{\mathbb{R}}$, $(f(t,\textbf{y},\mathbf{z}))_{t\in [0,T]}$ is
$\mathbf{F}$-progressively measurable for all
$(\mathbf{y},\mathbf{z})$ in ${\mathbb{R}}^{N+1} \times
{\mathbb{R}}^{(N+1)\times d}$. We make the following standard
assumptions, which extend naturally (A1) and (A2): \vskip0.2cm

(B1) There is some constant $C\ge 0$  such that, P-a.s., for all $t\in
[0, T],\ \mathbf{y}_{1}, \mathbf{y}_{2}\in {\mathbb{R}}^{N+1},\ \mathbf{z}_{1}, \mathbf{z}_{2}\in {\mathbb{R}}^{(N+1)\times d},\\
\mbox{ }\hskip1.2cm   |f(t, \mathbf{y}_{1},\mathbf{z}_{1}) -
f(t,\mathbf{y}_{2},\mathbf{z}_{2})|\leq
C(|\mathbf{y}_{1}-\mathbf{y}_{2}| +
|\mathbf{z}_{1}-\mathbf{z}_{2}|).$
 \vskip0.2cm

(B2) $f(\cdot,0,0)\in L^2_{\mathbf{F}}(0,T;{\mathbb{R}})$.

\bigskip

Let now $f:\Omega\times[0,T]\times {\mathbb{R}}^{N+1}\times
{\mathbb{R}}^{(N+1)\times d}\rightarrow {\mathbb{R}}$ be a measurable
function satisfying the assumptions (B1) and (B2). As above we
suppose that $L\in S^2_{\mathbf{F}}([0,T]),$\ $\xi\in
L^2(\Omega,{\cal F}_T,P)$\ and $L_T\leq\xi, \ P$-a.s. For a triplet
$(\xi,f,L)$ with these properties we say that it satisfies the
Standard Assumptions (B).

The above statements allow to extend the existence and uniqueness
result to RBSDEs whose data triplet $(\xi,f,L)$ satisfies the
Standard Assumptions (B).

\bp For every data triplet $(\xi,f,L)$
satisfying the Standard Assumptions (B) the RBSDE (N)
\be \begin{array}{ll}
& {\rm (i)}\ Y \in S_{\mathbf{F}}^2([0,T]), \, Z \in
L_{\mathbf{F}}^2([0,T];{\mathbb{R}}^d)\ \mbox{and}\ K\in
A_{\mathbf{F}}^{2,c}([0,T]);\\
& {\rm (ii)}\ Y_t = \displaystyle \xi +
\int_t^Tf(s,\Theta_N(Y_s,Z_s))ds + K_{T}-K_{t}
-\int^T_tZ_sdW_s,\quad t\in [0,T];\\
& {\rm (iii)}\ Y_t \geq L_t,\ \mbox{a.s., for all}\  t\in [0,T];\\
& {\rm (iv)}\ \displaystyle\int_0^T(Y_t - L_t)dK_{t}=0;
\end{array}
\ee
\noindent admits a unique solution $(Y,Z,K).$
\ep

\noindent\textbf{Proof}. Given an arbitrary couple $(U,V)\in
H^2=L_{\mathbf{F}}^2 ([0,T])\times
L_{\mathbf{F}}^2([0,T];{\mathbb{R}}^d)$ we put
$g(t)=f(t,\Theta_N(U_t,V_t)),\, t\in[0,T],$ and we denote by
$(Y,Z,K)$ the unique solution of the RBSDE with data triplet
$(\xi,g,L).$ For this we observe that the process $g$ is in
$L_{\mathbf{F}}^2([0,T])$ (and so it satisfies (A1) and (A2)) and we
recall that Lemma 3.1 guarantees the existence and the uniqueness of
the triplet $(Y,Z,K)$. We denote the mapping $(U,V)\rightarrow
(Y,Z)$ by $\Phi.$ For proving that the above RBSDE admits a unique
solution it suffices to show that, for a suitable equivalent norm in
$H^2$, the mapping $\Phi:H^2\rightarrow H^2$ is a contraction.
Indeed, if $\Phi$ is a contraction mapping on $H^2$ then there
exists a unique couple $(Y,Z)\in H^2$ such that $\Phi(Y,Z)=(Y,Z)$.
Due to the definition of $\Phi$, there is some $K\in
A_{\mathbf{F}}^{2,c}([0,T])$ such that $(Y,Z,K)$ is a solution of
the RBSDE with data triplet $(\xi,f(.,\Theta_N(Y,Z)),L).$
Consequently, $(Y,Z,K)$ is a solution of the above RBSDE. The
uniqueness of the solution of our RBSDE follows immediately from the
fact that whenever $(Y,Z,K)$ is a solution the couple $(Y,Z)$
is the unique fixed point of $\Phi$ in $H^2$.

For proving that the mapping $\Phi$ is a contraction with respect to
an appropriate equivalent norm on $H^2$ , we consider arbitrary $(U,V),(U',V')\in
H^2$ and apply Lemma 3.5 to $(Y,Z)=\Phi(U,V),$ $(Y',Z')=\Phi(U',V')$. For

\smallskip

$(\overline{Y},\overline{Z})=(Y-Y',Z-Z'),$
$(\overline{U},\overline{V})=(U-U',V-V')$, $\overline{\xi}=0,\,
\overline{L}=0,$

$\overline{g}(t)=f(t,\Theta_N(U_t,V_t))-f(t,\Theta_N(U'_t,V'_t)),\,
t\in[0,T],$

\smallskip

\noindent we thus get
$$
  \begin{array}{cl}
\displaystyle E\left[\int_0^Te^{\gamma t}(\vert\overline{Y}_t\vert^2+\vert \overline{Z}_t\vert^2) dt\right]\le & \displaystyle \delta E\left[ \int_0^Te^{\gamma t} \vert\overline{g}(t)\vert^2dt \right],\\
\end{array}
$$
for any $\delta>0$; the constant $\gamma>0$ depends only on $\delta$ and on the Lipschitz constant $L$ of $f(\omega,t,.,.)$. On the other hand,
$$
  \begin{array}{cl}
\displaystyle E\left[\vert\overline{g}(t)\vert^2\right] & \displaystyle \le L^2(N+1) \sum_{k=0}^NE \left[\vert\Theta^k (\overline{U}_t, \overline{V}_t) \vert^2 \right]\\
&
=\displaystyle L^2(N+1)^2E\left[\vert\overline{U}_t\vert^2+\vert\overline{V}_t\vert^2
\right], \quad t\in[0,T]
\end{array}
$$
(Recall that the random vectors $\Theta^k (\overline{U}_t,
\overline{V}_t),$ $k\ge 0,$ obey the same probability law).
Consequently, for $\delta:=\frac{1}{2} (L^2(N+1)^2)^{-1}$,
$$E\left[\int_0^Te^{\gamma t}(\vert\overline{Y}_t\vert^2+\vert \overline{Z}_t\vert^2) dt\right]\le \frac{1}{2}E\left[\int_0^Te^{\gamma t}(\vert\overline{U}_t\vert^2+\vert \overline{V}_t\vert^2) dt\right].$$
This shows that the mapping $\Phi:H^2\rightarrow H^2$ is contractive
with respect to the norm
$$\Vert(U,V)\Vert_{H^2}=\left(E[\int_0^T e^{\gamma t}(\vert U_t\vert^2+\vert V_t\vert^2) dt]\right)^{1/2},\, \,  (U,V)\in H^2.$$ The proof is complete.

\br Let us remark that for the type of
RBSDE which we have studied in Proposition 3.1 the comparison
principle does, in general, not hold. A consequence is that the
penalization method can't be used for the proof of the existence for
such a RBSDE.\er

Let us give an example:

\begin{example}
(1) We consider the BSDE without reflection
\be\label{3.8} Y_t=\xi+\int_t^T f(s, \Theta_N(Y_s, Z_s))ds-\int_t^TZ_sdW_s,\ \ t\in [0, T],\ee
with $\xi\in L^\infty (\Omega, {\cal{F}}_T^W, P)$,\ $f(s, y,z)=-y_1,\ y=(y_0, y_1,\cdots, y_N),\ z=(z_0, z_1, \cdots, z_N)$.\\
Then, the equation (\ref{3.8}) takes the form
$$ Y_t=\xi-\int_t^T  \Theta^1(Y_s)ds-\int_t^TZ_sdW_s,\ \ t\in [0, T].$$
Since $\Theta^1(Y_s),\ s\in [0, T]$, is independent of $\omega$, $Z$\ is obtained from the martingale representation property of $\xi\in L^\infty (\Omega, {\cal{F}}_T^W, P)$,
\be\label{3.9} \xi=E[\xi]+\int_0^T Z_sdW_s,\ \mbox{where}\ Z=(Z_s)\in L^2_{{\mathbb F}^W}(0, T);\ee
and $Y\in S^2_{{\mathbb F}}(0, T)$\ is the unique solution of the following equation:
\be\label{3.10} Y_t=E[\xi|{\cal{F}}_t^W]-\int_t^T  \Theta^1(Y_s)ds,\ \ t\in [0, T].\ee

Using the notation $\textbf{I}_i=(1,\cdots, 1)\in {\mathbb R}^i,\ i\geq 1$, and the fact that $\xi$\ as ${\cal{F}}_T^W$-measurable random variable coincides P-a.s. with some Borel measurable functional $\Phi: C_0([0, T])\rightarrow {\mathbb R}$\ combined with $W$, $\xi=\Phi(W)$, P-a.s., we see that the unique solution of (\ref{3.10}) is of the form
\be\label{3.11}\begin{array}{rcl}
  Y_t&=&E[\Phi(W)|{\cal{F}}_t^W]+\sum_{i=1}^\infty(-1)^i\int_t^T \int_{t_1}^T \cdots\int_{t_{i-1}}^TE[\Phi(W^{\textbf{I}_i})|{\cal{F}}_{t_i}^{W^{\textbf{I}_i}}]dt_i\cdots dt_2dt_1\\
&=& E[\Phi(W)|{\cal{F}}_t^W]+\sum_{i=1}^\infty(-1)^i\int_t^T\frac{(s-t)^{i-1}}{(i-1)!} E[\Phi(W^{\textbf{I}_i})|{\cal{F}}_s^{W^{\textbf{I}_i}}]ds, \ t\in [0, T].
\end{array}\ee
Indeed, due to the definition of $\Theta^k, \ k\geq 1,$
$$
  \Theta^1(Y_s)=E[\Phi(W^{(1)})|{\cal{F}}_s^{W^{(1)}}]+\sum_{i=1}^\infty(-1)^i\int_s^T \int_{t_1}^T \cdots\int_{t_{i-1}}^TE[\Phi(W^{\textbf{I}_{i+1}})|{\cal{F}}_{t_i}^{W^{\textbf{I}_{i+1}}}]dt_i\cdots dt_2dt_1, \ s\in [0, T],\\
$$
and it can be easily checked that $(Y, \Theta^1(Y))$\ satisfies (\ref{3.10}), and (\ref{3.10}) with (\ref{3.9}) yields (\ref{3.8}). Consequently, $Y$\ given by (\ref{3.11}) and $Z$\ by (\ref{3.9})\ is the unique solution of (\ref{3.8}). We also observe that, if $|\xi|\leq C$, P-a.s., then
$$|Y_t|\leq C\sum_{i=0}^\infty\frac{(T-t)^i}{i!}=Ce^{T-t},\ t\in[0, T].$$
Consequently, (\ref{3.8}) can be regarded also as an RBSDE with reflection barrier $L_t=-Ce^T,\ t\in [0, T]$, and its unique solution $(Y, Z, K)$ is given by (\ref{3.11}), (\ref{3.10}) and $K_t=0$, $t\in [0, T]$.

(2) Let us now consider the RBSDE introduced above with $T=2$, and $\xi=|W_1|^2\wedge 1$, $C=1,\ L_t=-e^T,\ t\in [0, T]$. Then, again $K_t=0,\ t\in [0, T]$, and from (\ref{3.8}),
$$\begin{array}{rcl}
  E[Y_t]&=&E[|W_1|^2\wedge 1]-\int_t^2E[\Theta^1(Y_s)]ds\\
&=& E[|W_1|^2\wedge 1]-\int_t^2E[Y_s]ds, \ t\in [0, 2],\\
\end{array}
$$
i.e., $E[Y_t]=E[|W_1|^2\wedge 1]e^{-(2-t)},\ t\in [0, 2].$\ On the other hand, since $\Theta^1(Y_s),\ s\in [0, 2]$, is independent of $W$,
$$ E[Y_t|{\cal{F}}_T^{W}]=|W_1|^2\wedge 1-\int_t^2E[Y_s]ds,\ t\in [0, 2],$$
and thus,
\be\label{3.5}\begin{array}{rcl}
E[Y_1|{\cal{F}}_T^{W}]&=&|W_1|^2\wedge 1-\int_1^2E[Y_s]ds\\
&=& |W_1|^2\wedge 1-E[|W_1|^2\wedge 1](1-e^{-1})\\
&<& 0,\ \ \mbox{on}\ \{|W_1|^2\wedge 1<E[|W_1|^2\wedge 1](1-e^{-1})\}.
\end{array}
\ee
Consequently, $P\{E[Y_1|{\cal{F}}_T^{W}]<0\}>0$, and, hence, also $P\{Y_1<0\}>0$.

On the other hand, for the terminal condition $\xi'=0$, $(Y', Z', K')=(0, 0, 0)$\ is the unique solution of our RBSDE. This shows that, although $P\{\xi>\xi'\}=1$, we have $P\{Y_1<Y'_1\}>0$, i.e., in general, our RBSDE doesn't satisfy a comparison principle.
\end{example}

\section{A Limit Approach for Mean-Field BSDEs with Reflection}

The objective of this section is to study the limit of RBSDE(N) as
$N$ tends to infinity. For this we choose the framework we have
already introduced for the study of the approximation of the reflected mean-field BSDE
by RBSDEs. Let $(\Phi,g,{\cal X})$ be a data triplet satisfying the
assumptions (C1)-(C3):
\smallskip

(C1) $g:\Omega\times [0,T]\times\left({\mathbb{R}}^m\times
{\mathbb{R}}\times {\mathbb{R}}^d\right)^2\rightarrow {\mathbb{R}}$
is a bounded measurable function, and $g$ is Lipschitz with respect to
$(\mathbf{u,v})$, i.e., $P$-a.s., for
all $t\in[0,T]$ and
$(\mathbf{u,v}),(\mathbf{u',v'})\in\left({\mathbb{R}}^m\times
{\mathbb{R}}\times{\mathbb{R}}^d\right)^2$,
$$\vert g(t,(\mathbf{u,v}))-g(t,(\mathbf{u',v'}))\vert\le C\left(\vert \mathbf{u}-
\mathbf{u'}\vert+\vert \mathbf{v}-\mathbf{v'}\vert\right);$$

\smallskip

(C2) $\Phi:\Omega\times
{\mathbb{R}}^m\times{\mathbb{R}}^m\rightarrow {\mathbb{R}}$ is a
bounded measurable function, and $\Phi(\omega,.,.)$ is
Lipschitz, i.e., $P$-a.s., for all
$(x,\hat{x}), (x',\hat{x}')\in\mathbb{R}^m,$
$$\vert\Phi(x,\hat{x})-\Phi(x',\hat{x}')\vert\le C\left(\vert x-x'\vert+\vert\hat{x}-\hat{x}'\vert \right).$$

\smallskip

(C3) ${\cal X}=(X^N)_{N\ge 1}$ is a Cauchy sequence in ${\cal
S}_{\mathbf{F}}^2([0,T];{\mathbb{R}}^m)$, i.e., there is a (unique)
process $X\in {\cal S}_{\mathbf{F}}^2([0,T];{\mathbb{R}}^m)$ such
that
$$\displaystyle E[\sup_{t\in[0,T]}\vert X^N_t-X_t\vert^2]\rightarrow
0\ \mbox{as}\ N\rightarrow +\infty.$$ \vskip0.5cm

\noindent Moreover, let $h:\Omega\times [0,T]\times {\mathbb{R}}^m\rightarrow
{\mathbb{R}}$ be a function with the following properties: \smallskip

(C4) $h:\Omega\times [0,T]\times {\mathbb{R}}^m\rightarrow
{\mathbb{R}}$ is a bounded measurable function which is

$\bullet$ $\mathbf{F}$-progressively measurable, for every fixed
$x\in {\mathbb{R}}^m$;

$\bullet$ Lipschitz continuous, for every fixed
$(\omega,t)\in\Omega\times[0,T],$ with a Lipschitz constant that
doesn't depend on $(\omega,t)$;

$\bullet$ continuous in $t$, for every fixed
$(\omega,x)\in\Omega\times {\mathbb{R}}^m$.
 \vskip0.2cm

Given such a quadruplet $(\Phi,g,h,{\cal X})$ such that
$(\Phi,g,{\cal X})$  fulfills the assumptions (C1)-(C3), $h$\
satisfies the assumption (C4) and $h(T,x)\le\Phi(x,x'),$
$P$-a.s., for all $(x,x')\in {\mathbb{R}}^m\times {\mathbb{R}}^m$,
we say that $(\Phi,g,h,{\cal X})$ satisfies the Standard Assumptions
(C) and we put, for $N\ge 1$ and $\omega\in\Omega$,
$$
  \begin{array}{cl}
& \displaystyle \xi^N(\omega):=\frac{1}{N}\sum_{k=1}^N\Phi(\Theta^k(\omega),X^N_T(\omega),X^N_T(\Theta^k(\omega))),\\
& \displaystyle f^N(\omega,t,\textbf{y},\textbf{z}):= \frac{1}{N}\sum_{k=1}^N g(\Theta^k(\omega),t,X^N_t(\omega),(y_0,z_0),X^N_t(\Theta^k(\omega)),(y_k,z_k)),\\
&\mbox{\quad for }t\in[0,T],\, \textbf{y}=(y_0,\cdots,y_N)\in {\mathbb{R}}^{N+1},\textbf{z}=(z_0,\cdots,z_N)\in {\mathbb{R}}^{(N+1)\times d},\\
& \\
& \displaystyle L^N_t(\omega):=h(\omega,t,X^N_t(\omega)),\,
t\in[0,T].
\end{array}
$$

\medskip

We notice that, for each $N\ge 1$, the triplet $(\xi^N,f^N,L^N)$
satisfies the Standard Assumptions for an RBSDE: $(\xi^N,f^N)$ satisfies (B1)-(B2) and $L^N\in {\cal S}^2_{\mathbf{F}}([0,T])$ is such that $L^N_T\le \xi^N$. Thus, due to Proposition
3.1, we have for all $N\ge 1$\ a unique solution
$(Y^N,Z^N,K^N)$ of the RBSDE (N)
\be
  \begin{array}{ll}
& {\rm (i)}\ Y^N \in S_{\mathbf{F}}^2([0,T]), \, Z^N \in
L_{\mathbf{F}}^2([0,T];{\mathbb{R}}^d)\ \mbox{and}\ K^N\in
A_{\mathbf{F}}^{2,c}([0,T]);\\
& {\rm (ii)}\ Y^N_t = \displaystyle \xi^N +
\int_t^Tf^N(s,\Theta_N(Y_s,Z_s))ds + K^N_{T}-K^N_{t}
-\int^T_tZ^N_sdW_s,\ t\in [0,T];\\
& {\rm (iii)}\ Y^N_t \geq L^N_t,\ \mbox{a.s., for any}\  t\in [0,T];\\
& {\rm (iv)}\ \displaystyle\int_0^T(Y^N_t - L^N_t)dK^N_{t}=0.
\end{array}
\ee

\noindent We remark that, the driving coefficient of the above
RBSDE(N) can be written as follows:
$$f^N(s,\Theta_N(Y_s^N,Z_s^N))=\frac{1}{N}\sum_{k=1}^N(\Theta^kg)(s,(X_s^N,Y_s^N, Z_s^N),\Theta^k(X_s^N,Y_s^N,Z_s^N)),$$
$s\in[0,T].$ Our objective is to show that the unique solution of
RBSDE(N) converges to the unique solution $(Y,Z,K)$ of the Reflected
Mean-Field BSDE
\be
  \begin{array}{ll}
&{\rm (i)}\ Y \in S_{\mathbf{F}}^2([0,T]), \, Z \in
L_{\mathbf{F}}^2([0,T];{\mathbb{R}}^d)\ \mbox{and}\  K\in
A_{\mathbf{F}}^{2,c}([0,T]);\\
& {\rm (ii)}\ Y_t = \displaystyle E\left[\Phi(x,X_T)\right]_{x=X_T} +
\int_t^TE\left[g(s,\mathbf{u},\Lambda_s)\right]_{\vert
\mathbf{u}=\Lambda_s}ds + K_{T}-K_{t} -\int^T_tZ_sdW_s;\\
& {\rm (iii)}\ \displaystyle Y_t \geq h(t,X_t),\ \mbox{a.s., for all}\  t\in
[0,T];\\
& {\rm (iv)}\ \displaystyle\int_0^T(Y_t - h(t,X_t))dK_{t}=0,\\
\end{array}
\ee

\noindent where we have used the notation $\Lambda=(X,Y,Z).$

\medskip

\bl Under the Standard Assumptions (C) on
the data quadruplet $(\Phi,g,h,{\cal X})$ the above  Reflected  Mean-Field BSDE
possesses a unique solution $(Y,Z,K)\in
S_{\mathbf{F}}^2([0,T])\times
L_{\mathbf{F}}^2([0,T];{\mathbb{R}}^d)\times
A_{\mathbf{F}}^{2,c}([0,T]).$\el
\noindent The proof is standard. For the convenience we give the proof here.

\medskip

\noindent\textbf{Proof}. Let
$H^2:=L_{\mathbf{F}}^2([0,T];{\mathbb{R}})\times
L_{\mathbf{F}}^2([0,T];{\mathbb{R}}^d).$\ Similar to the discussion in the beginning of the proof of Proposition 3.1 it is sufficient to prove
the existence and the uniqueness for the above BSDE in $H^2$.
Indeed, if $(Y,Z)$ is a solution of our BSDE in $H^2$, an easy
standard argument shows that it is also in
$B^2:=S_{\mathbf{F}}^2([0,T];{\mathbb{R}})\times
L_{\mathbf{F}}^2([0,T];{\mathbb{R}}^d).$ On the other hand, the
uniqueness in $H^2$ implies obviously that in its subspace $B^2.$

For proving the existence and uniqueness in $H^2$ we consider for an
arbitrarily given couple of processes $(U,V)\in H^2$ the coefficient
$g_s^{U,V}=E[g(s,\lambda, \Lambda_s)]|_{\lambda=\Lambda_s},\, s\in[0,T],$ for $\Lambda_s=(X_s, U_s, V_s)$. Since $g$ is an element
of $L_{\mathbf{F}}^2([0,T];{\mathbb{R}})$ it follows from Lemma 3.1
that there is a unique solution $\Phi(U,V):=(Y,Z)\in H^2$ of the reflected BSDE:
$$
  \begin{array}{ll}
&{\rm (i)}\ Y \in S_{\mathbf{F}}^2([0,T]), \, Z \in
L_{\mathbf{F}}^2([0,T];{\mathbb{R}}^d)\ \mbox{and}\  K\in
A_{\mathbf{F}}^{2,c}([0,T]);\\
& {\rm (ii)}\ Y_t = \displaystyle E\left[\Phi(x,X_T)\right]_{x=X_T} +
\int_t^Tg_s^{U,V}ds + K_{T}-K_{t} -\int^T_tZ_sdW_s;\\
& {\rm (iii)}\ \displaystyle Y_t \geq h(t,X_t),\ \mbox{a.s., for all}\  t\in
[0,T];\\
& {\rm (iv)}\ \displaystyle\int_0^T(Y_t - h(t,X_t))dK_{t}=0,\\
\end{array}
$$
For a such defined mapping $\Phi:H^2\rightarrow H^2$ it suffices to
prove that it is a contraction with respect to an appropriate
equivalent norm on $H^2$, in order to complete the proof. For this
end, we consider two couples $(U^1,V^1),(U^2,V^2)\in H$ and
$(Y^k,Z^k)=\Phi(U^k,V^k),\, k=1,2.$ Then, due to Lemma 3.5, for all
$\delta>0$ there is some constant $\gamma>0$ (only depending on
$\delta$) such that, with the notation
$(\overline{Y},\overline{Z})=(Y^1-Y^2,Z^1-Z^2),$
$$\begin{array}{cl}
  & \displaystyle E[\int_0^T e^{\gamma t}(\vert\overline{Y}_t\vert^2+\vert\overline{Z}_t \vert^2)dt] \\
  & \le\displaystyle \delta E[\int_0^T e^{\gamma t}\vert g_t^{U^1,V^1}- g_t^{U^2,V^2}\vert^2dt].
  \end{array}$$
Let $(\overline{U},\overline{V})=(U^1-U^2,V^1-V^2).$ Then, from the
Lipschitz continuity (C1) of $g$ (with Lipschitz constant
$C$ which doesn't depend on $(\omega,t)$)
$$\begin{array}{cl}
   & \displaystyle E[|g_t^{U^1,V^1}-g_t^{U^2,V^2}|^2]\\
  & =\displaystyle E[|E[g(t, \lambda^1, \Lambda^1_t)]|_{\lambda^1=\Lambda^1_t}-E[g(t, \lambda^2, \Lambda^2_t)]|_{\lambda^2=\Lambda^2_t}|^2]\\
  & \displaystyle\le C E[|\Lambda^1_t-\Lambda^2_t|^2] \\
  &=CE[|U_t^1-U_t^2|^2+|V_t^1-V_t^2|^2],
  \end{array}$$
where $\Lambda^1_t:=(X_t,U^1_t,V^1_t),\ \Lambda^2_t:=(X_t,U^2_t,V^2_t).$
Consequently, we have
$$\begin{array}{cl}
  & \displaystyle E[\int_0^T e^{\gamma t}(\vert\overline{Y}_t\vert^2+\vert\overline{Z}_t \vert^2)dt] \le\displaystyle \delta C \int_0^T e^{\gamma t}E[\vert \overline{U}_t\vert^2+\vert\overline{V}_t\vert^2]dt\\
  &=\displaystyle\frac{1}{2}E[\int_0^T e^{\gamma t} (\vert \overline{U}_t\vert^2+\vert\overline{V}_t\vert^2) dt],
  \end{array}$$
for $\delta:=\frac{1}{2C}$. This shows that if we
endow the space $H^2$ with the norm
$$\Vert(U,V)\Vert_{H^2}=\left(E[\int_0^T e^{\gamma t}(\vert U_t\vert^2+\vert V_t\vert^2) dt]\right)^{1/2},\, \,  (U,V)\in H^2,$$ the mapping $\Phi:H^2\rightarrow H^2$ becomes a contraction. Thus, the proof is complete.\endpf

\smallskip

We now can formulate the following theorem:
\bt Under the Standard Assumptions (C) on
the data quadruplet $(\Phi,g,h,{\cal X})$, the unique solution
$(Y^N,Z^N,K^N)$ of RBSDE(N) (4.1) converges to the unique solution
$(Y,Z,K)$ of the above MFBSDE (4.2) with reflection:
$$E\left[\sup_{t\in[0,T]}\vert Y^N_t-Y_t\vert^2+\int_0^T\vert Z^N_t-Z_t\vert^2dt+\sup_{t\in[0,T]}\vert K^N_t-K_t\vert^2\right]\rightarrow 0,$$
$\mbox{ as } N\rightarrow +\infty.$\et

\noindent\textbf{Proof}. First we want to prove that

\noindent\underline{Step 1}.  For all $p\ge 2,$

$$E\left[\int_0^T\left\vert\frac{1}{N}\sum_{k=1}^N(\Theta^kg)(t,\Lambda_t, \Theta^k(\Lambda_t)) - E\left[g(t,\mathbf{u},\Lambda_t)\right]_{\vert \mathbf{u}=\Lambda_t}\right\vert^pdt\right]\rightarrow 0,$$
and
$$E\left[\left\vert\frac{1}{N}\sum_{k=1}^N(\Theta^k\Phi)(X_T,\Theta^k(X_T))- E\left[\Phi(x,X_T)\right]\vert_{x=X_T}\right\vert^p\right]\rightarrow 0,$$
as $N\rightarrow +\infty$ (notice that
$(\Theta^k\Phi)(\omega,X_T,\Theta^k(X_T))(\omega):=\Phi(\Theta^k(\omega),
X_T(\omega),X_T(\Theta^k(\omega)))$).

To prove the first convergence we need to consider arbitrary $t\in[0,T]$
and $\mathbf{u}\in {\mathbb{R}}^m\times {\mathbb{R}}\times
{\mathbb{R}}^d.$ Notice that the sequence of random variables
$(\Theta^kg)(t,\mathbf{u},\Lambda_t),\ k\ge 1,$ is i.i.d. and has the
same law as $g(t,\mathbf{u},\Lambda_t),$  from the Strong
Law of Large Numbers we get that
$$\frac{1}{N}\sum_{k=1}^N (\Theta^kg)(t,\mathbf{u},\Lambda_t)\longrightarrow E\left[ g(t,\mathbf{u},\Lambda_t)\right],$$
$P$-a.s., as $N\rightarrow +\infty.$ For an arbitrarily
small $\varepsilon>0,$ let $\Lambda_t^\varepsilon:\Omega\rightarrow
{\mathbb{R}}^m\times {\mathbb{R}}\times {\mathbb{R}}^d$ be a random
vector which has only a countable number of values, and also satisfies
$\vert\Lambda_t-\Lambda_t^\varepsilon\vert\le\varepsilon,$
everywhere on $\Omega$. Then, obviously, $$\frac{1}{N}\sum_{k=1}^N
(\Theta^kg)(t,\Lambda_t^\varepsilon,\Theta^k(\Lambda_t^\varepsilon))\longrightarrow
E\left[
g(t,\mathbf{u},\Lambda_t^\varepsilon)\right]_{\vert\mathbf{u}=\Lambda_t^\varepsilon},$$
\noindent $P$-a.s., as $N$ tends to $+\infty.$ On the other hand, from the
Lipschitz continuity of $g(\omega,t,.,\mathbf{v}),$ uniformly in $(\omega,t,\mathbf{v})$, we know that $\Lambda_t$ also have the
convergence:
$$\frac{1}{N}\sum_{k=1}^N(\Theta^kg)(t,\Lambda_t,\Theta^k(\Lambda_t)) \longrightarrow E\left[g(t,\mathbf{u},\Lambda_t)\right]_{\vert \mathbf{u}=\Lambda_t},$$
$P$-a.s., as $N\rightarrow +\infty.$ Finally, from the
boundedness of $g$ and, thus, of that of the convergence, we get
the wished result. Similarly, we also obtain the
$L^p$-convergence for the terminal conditions, for all $p\ge 2.$

\medskip

\noindent\underline{Step 2}. Recalling the argument given in Step 1
we see that, for all $p\ge 2,$ as $N\rightarrow +\infty,$
$$\begin{array}{cl}
  & \displaystyle {{R}}^{N,p}_1:= E\left[\left\vert\frac{1}{N}\sum_{k=1}^N(\Theta^k\Phi) (X_T, \Theta^k(X_T))- E\left[\Phi(x,X_T)\right]_{x=X_T}\right\vert^p\right] \rightarrow 0,\\
  & \mbox{ and } \\
  & \displaystyle {{R}}^{N,p}_2:=E\left[\int_0^T\left\vert\frac{1}{N}\sum_{k=1}^N (\Theta^kg)(t,\Lambda_t, \Theta^k(\Lambda_t)) - E\left[g(t,\mathbf{u}, \Lambda_t)\right]_{\vert \mathbf{u}=\Lambda_t}\right\vert^pdt\right]\rightarrow 0.
  \end{array}$$
For estimating the distance between $(Y^N,Z^N,K^N)$ and $(Y,Z,K)$ we
apply Lemma 3.5. and get, for $\delta\in(0,1)$  which will be
specified later, and for some $\gamma>0$ (depending on $\delta$ and
on the Lipschitz constant of $g$),
$$
  \begin{array}{cl}
& \displaystyle E\left[\int_0^Te^{\gamma t}(\vert Y^N_t-Y_t\vert^2+\vert Z^N_t -Z_t\vert^2)dt\right]\le
2 \displaystyle E[e^{\gamma T}\vert\overline{\xi}^N\vert^2] +\\
& \mbox{ }\hskip0.3cm \displaystyle
+\delta E\left[\int_0^Te^{\gamma t}\vert \overline{f}^N(t,Y_t,Z_t)\vert^2 dt\right]
+\left(C E\left[\sup_{t\in[0,T]}\vert e^{\gamma t}\overline{L}^N_t\vert^2\right]\right)^{1/2} (\Psi_{0,T}^N)^{1/2},\\
  \end{array}
$$
where
$$
  \begin{array}{rcl}
\displaystyle \overline{\xi}^N& := & \displaystyle\frac{1}{N}\sum_{k=1}^N(\Theta^k \Phi)(X^N_T, \Theta^k(X^N_T))-E\left[\Phi(x,X_T)\right]_{x=X_T},\\
\displaystyle \overline{f}^N(t,Y_t,Z_t) & := & \displaystyle\frac{1}{N}\sum_{k=1}^N (\Theta^kg)(t,\Lambda^N_t,\Theta^k(\Lambda^N_t))-E\left[g(t,\mathbf{u},\Lambda_t) \right]_{\vert \mathbf{u}=\Lambda_t},\\
& &\mbox{ }\hskip 0.5cm \hbox{with } \Lambda^N_t=(X^N_t,Y^N_t,Z^N_t),\, \Lambda_t=(X_t, Y_t,Z_t),\\
\displaystyle \overline{L}^N_t & := & h(t,X_t^N)-h(t,X_t),
  \end{array}
$$
and
$$
  \begin{array}{cl}
& \displaystyle\Psi^N_{t,T} := \displaystyle E\bigg[\left\vert\xi^N\right\vert^2+ \left\vert E\left[\Phi(x,X_T)\right]_{x=X_T}\right\vert^2 +\left(\int_t^T\left\vert f^N(s,\Theta_N(Y^N_s,Z^N_s))\right\vert ds\right)^2\\
& \displaystyle +\left(\int_t^T\left\vert E\left[g(s,\mathbf{u},\Lambda_t) \right]_{\vert \mathbf{u}= \Lambda_s}\right\vert ds\right)^2+ \sup_{s\in[t,T]}\vert h(s,X^N_s)\vert^2+\sup_{s\in[t,T]}\vert h(s,X_s)\vert^2\vert{\cal F}_t\bigg],\\
  \end{array}
$$
for $t\in[0,T].$ In virtue of the boundedness of the coefficients
$\Phi,\, g$ and $h$ it follows that, for some constant $C$,
$\Psi^N_{t,T}\le C, \, t\in[0,T],$ $P$-a.s. On the other hand,
recalling that the coefficients $\Phi(\omega,.)$, $g(\omega,t,.,.)$
and $h(\omega,t,.)$ are Lipschitz, with some Lipschitz constant $L$
which is independent of $(\omega,t)$, and using the the fact that, for any
random variable $\xi$, the variables $\Theta^k(\xi),$ $k\ge 0$, obey
the same probability law, we see that
$$
  \begin{array}{rcl}
\displaystyle E\left\vert\overline{\xi}^N\right\vert^2\vert] & \le & 2R^{N,2}_1+ 8L^2 E\left[\left\vert X^N_T-X_T\right\vert^2\right],\\
\displaystyle E\left[\int_0^Te^{\gamma t}\vert \overline{f}^N(t,Y_t,Z_t)\vert^2 dt\right] & \le & \displaystyle 2e^{\gamma T}R^{N,2}_2+8L^2E\left[\int_0^Te^{\gamma t}\vert\Lambda^N_t- \Lambda_t\vert^2dt\right],\\
 \displaystyle E\big[\sup_{t\in[0,T]}\vert \overline{L}^N_t\vert^2\big] & \le & L^2 E\big[\sup_{t\in[0,T]}\vert X^N_t-X_t\vert^2\big].\\
  \end{array}
$$
Consequently, with the notation
$$
  \begin{array}{cl}
R_N & := \displaystyle 4e^{\gamma T}\bigg( R_1^{N,2}+R_2^{N,2}+4L^2E[\vert X_T^N-X_T \vert^2]\\
& \displaystyle +2L^2E\left[\int_0^T\vert X_t^N-X_t\vert^2dt\right]+ CL\left(E\left[ \sup_{t\in[0,T]}\vert X_t^N-X_t\vert^2\right]\right)^{1/2}\bigg)\\
  \end{array}
$$
we have
$$
  \begin{array}{cl}
& \displaystyle E\left[\int_0^Te^{\gamma t}(\vert Y^N_t-Y_t\vert^2+\vert Z^N_t -Z_t\vert^2)dt\right]\\
& \le \displaystyle R_N+8L^2\delta \displaystyle E\left[\int_0^Te^{\gamma t}(\vert Y^N_t-Y_t\vert^2+\vert Z^N_t -Z_t\vert^2)dt\right],\\
  \end{array}
$$
and choosing $\delta:=\frac{1}{16}L^{-2}$ we obtain
$$  E\left[\int_0^Te^{\gamma t}(\vert Y^N_t-Y_t\vert^2+\vert Z^N_t -Z_t\vert^2)dt\right]\le 2R_N,\, N\ge 1.$$
Hence, since $R_N$ converges to zero as $N$ tends to $+\infty,$ we
also have
$$ E\left[\int_0^Te^{\gamma t}(\vert Y^N_t-Y_t\vert^2+\vert Z^N_t -Z_t\vert^2)dt\right]\longrightarrow 0,\mbox{ as } N\rightarrow +\infty.$$
Applying now Lemma 3.4 we obtain the following estimate, for all
$t\in[0,T],$ $P$-a.s.,
$$
  \begin{array}{cl}
& \displaystyle E[\sup_{s\in[t,T]}\vert Y^N_s-Y_s\vert^2+\int_t^T\vert Z^N_s -Z_s\vert^2ds+\vert (K^N_T-K_T)-(K^N_t-K_t)\vert^2]\\
& \le C \displaystyle E[\vert\overline{\xi}^N\vert^2 +(\int_t^T\vert \overline{f}^N(s,Y_s,Z_s)\vert ds)^2]+C E[\sup_{s\in[t,T]}\vert\overline{L}^N_s\vert^2]^{1/2} E[\Psi_{t,T}^N]^{1/2},\\
  \end{array}
$$
which right-hand side converges to zero according to our preceding
convergence result. Consequently, $E[\sup_{s\in[t,T]}\vert
Y^N_s-Y_s\vert^2]\rightarrow 0$ and, for all $t\in[0,T],$ $E[\vert
K_t^N-K_t\vert^2]\rightarrow 0$, as $N$ tends towards $+\infty$. In order to
conclude, it suffices to observe that the fact that $K$ is a square
integrable, increasing continuous process implies that we even have
$$E[\sup_{t\in[0,T]}\vert K_t^N-K_t\vert^2]\rightarrow 0,\mbox{ as }N\rightarrow 0.$$ Indeed, given an arbitrary $\varepsilon>0$ we can find some finite partition $0=t_0<t_1<\cdots<t_M=T$ such that $E[\max_{1\le i\le M}(K_{t_i}-K_{t_{i-1}})^2]\le\varepsilon^2.$ Then, since the processes $K^N$ and $K$ are increasing,
$$
  \begin{array}{rcl}
\displaystyle \sup_{t\in[0,T]}\vert K^N_t-K_t\vert & = & \displaystyle \max_{1\le i\le M}\left(\sup_{t\in[t_{i-1},t_i]}\vert K^N_t-K_t\vert\right)\\
& \le & \displaystyle \max_{1\le i\le M}\left(\vert K^N_{t_{i-1}}-K_{t_i}\vert + \vert K^N_{t_i}-K_{t_{i-1}}\vert\right),\\
  \end{array}
$$
and, consequently,
$$
 \begin{array}{cl}
& \displaystyle E\left[\sup_{t\in[0,T]}\vert K^N_t-K_t\vert^2\right]
 \le  \displaystyle 2\sum_{1\le i\le M}\left(E[\vert K^N_{t_{i-1}}-K_{t_i}\vert^2] + E[\vert K^N_{t_i}-K_{t_{i-1}}\vert^2]\right)\\
& \mbox{ }\hskip0.2cm \rightarrow \displaystyle 4\sum_{1\le i\le M}E[( K_{t_i}-K_{t_{i-1}})^2]\le 4E[K_T\max_{1\le i\le M}(K_{t_i}-K_{t_{i-1}})]\\
& \mbox{ }\hskip0.5cm\le 4\left(E[K_T^2]\right)^{1/2}\varepsilon,
\mbox{ as } N\rightarrow +\infty.
\end{array}
$$
The proof is complete.

\section{Existence of a solution of the Reflected MFBSDE: approximation via penalization}

In~\cite{EKPPQ} the penalization method for BSDEs is used to prove the existence for the reflected BSDE. Can we use also here this method adapted to mean-field BSDEs, in order to study reflected MFBSDEs? In this section we will give a positive answer to this question. We will see that we can get the reflected mean-field BSDEs by using the penalization method to the mean-field BSDEs. The result of this section will be very useful in
Section 6.

For the given coefficient $g(\omega, t, y, z, \tilde{y}):  {\Omega}\times[0,T]\times
{\mathbb{R}} \times {\mathbb{R}}^{d}\times {\mathbb{R}}  \rightarrow {\mathbb{R}}$ which satisfies (A3), and $g$\ is nondecreasing with respect to $\tilde{y}$, the obstacle process $L\in S_{\mathbf{F}}^2([0,T])$, and the terminal condition $\xi\in L^2(\Omega, {\cal F}_T, \mathbb{R})$ such that $\xi\geq L_T$, P-a.s., we consider the following reflected mean-field BSDE:
\be\label{5.1}
\begin{array}{rcl}
& & {\rm (i)} Y \in S_{\mathbf{F}}^2([0,T]), \, Z \in
L_{\mathbf{F}}^2([0,T];{\mathbb{R}}^d)\ \mbox{and}\ K\in
A_{\mathbf{F}}^{2,c}([0,T]);\\
& & {\rm (ii)} Y_t = \displaystyle \xi + \int_t^TE\left[g(s,y,z,Y_s)\right]_{\vert
y=Y_s, z=Z_s}ds + K_{T}-K_{t} -\int^T_t Z_sdW_s;\\
& & {\rm (iii)} \displaystyle Y_t \geq L_t,\ \mbox{a.s.,\ for all}\ t\in
[0,T]; \\
& & {\rm (iv)} \displaystyle\int_0^T(Y_t - L_t)dK_{t}=0.
\end{array}
\ee

Similarly to the proof of Lemma 4.1 we know that the above equation (\ref{5.1}) has a unique solution $(Y, Z, K)$.

We define $f(s, y, z):= E[g(s, y, z, Y_s)]$, where $(Y, Z, K)$\ is the solution of (\ref{5.1}), $s\in [0, T], y\in {\mathbb R}, z\in {\mathbb R}^d.$

For each $n\in {\mathbb N}$, let $({Y}^n, {Z}^n)\in S_{\mathbf{F}}^2([0,T])\times L_{\mathbf{F}}^2([0,T];{\mathbb{R}}^d)$\ denote the solution of the following BSDE which in fact, also can be seen as a special mean-field BSDE:

\be\label{5.2} {Y}_t^n=\xi+\int_t^T f(s, {Y}_s^n, {Z}_s^n)ds+n\int_t^T({Y}_s^n-L_s)^{-}ds-\int_t^T{Z}_s^ndW_s,\ee
Then from the comparison theorem--Lemma 2.3 we know

\be\label{5.3} Y_t^n\leq Y_t^{n+1},\ 0\leq t\leq T,\ \mbox{a.s.},\  \forall n\in {\mathbb N}.\ee

\noindent We define $K_t^n:=n\int_0^t({Y}_s^n-L_s)^{-}ds,\ 0\leq t\leq T$.

\noindent Therefore, from the proof on Pages 719-723 in ~\cite{EKPPQ}, we know:
\be\label{5.4}\begin{array}{lll}
&{\rm(i)}\ Y_t^n\uparrow \widetilde{Y}_t, 0\leq t\leq T, \mbox{a.s.};\\
&{\rm(ii)}\ E\(\sup_{0\leq t\leq T}|Y_t^n-\widetilde{Y}_t|^2+\int_0^T|Z_s^n-\widetilde{Z}_s|^2ds+\sup_{0\leq t\leq T}|K_t^n-\widetilde{K}_t|^2\)\rightarrow 0,
\end{array}
\ee

\noindent where $(\widetilde{Y}, \widetilde{Z}, \widetilde{K})$\ is the solution of the following reflected BSDE:
\be\label{5.5}
\begin{array}{rcl}
& & {\rm (i)} \widetilde{Y} \in S_{\mathbf{F}}^2([0,T]), \, \widetilde{Z} \in
L_{\mathbf{F}}^2([0,T];{\mathbb{R}}^d)\ \mbox{and}\ \widetilde{K}\in
A_{\mathbf{F}}^{2,c}([0,T]);\\
& & {\rm (ii)} \widetilde{Y}_t = \displaystyle \xi + \int_t^Tf(s,\widetilde{Y}_s,\widetilde{Z}_s)ds + \widetilde{K}_{T}-\widetilde{K}_{t} -\int^T_t \widetilde{Z}_sdW_s;\\
& & {\rm (iii)} \displaystyle \widetilde{Y}_t \geq L_t,\ \mbox{a.s.,\ for all}\ t\in
[0,T]; \\
& & {\rm (iv)} \displaystyle\int_0^T(\widetilde{Y}_t - L_t)d\widetilde{K}_{t}=0.
\end{array}
\ee

By comparing the equations (\ref{5.1}) and (\ref{5.5}) we get from the uniqueness of the solution of RBSDE (\ref{5.5}) that

\be\label{5.6} Y_t=\widetilde{Y}_t,\ K_t=\widetilde{K}_t,\ 0\leq t\leq T,\ \mbox{a.s.},\ \mbox{and}\ \ Z_t=\widetilde{Z}_t,\ \mbox{a.e.a.s.}\ee

\noindent Therefore, we know:
\be\label{5.7}\begin{array}{lll}
&{\rm(i)}\ Y_t^n\uparrow {Y}_t, 0\leq t\leq T, \mbox{a.s.};\\
&{\rm(ii)}\ E\(\sup_{0\leq t\leq T}|Y_t^n-{Y}_t|^2+\int_0^T|Z_s^n-{Z}_s|^2ds+\sup_{0\leq t\leq T}|K_t^n-{K}_t|^2\)\rightarrow 0,
\end{array}
\ee
where $(Y, Z, K)$\ is the solution of reflected MFBSDE (\ref{5.1}).

Now we consider the following penalized mean-field BSDEs:
\be\label{5.8} \widetilde{{Y}}_t^n=\xi+\int_t^T E[g(s, y, z, \widetilde{Y}^n_s)]|_{y=\widetilde{Y}^n_s, z=\widetilde{Z}^n_s}ds+n\int_t^T(\widetilde{{Y}}_s^n-L_s)^{-}ds-\int_t^T\widetilde{{Z}}_s^ndW_s.\ee

\noindent From Lemmas 2.4 and 2.5 we know that it has a unique solution $(\widetilde{Y}^n, \widetilde{Z}^n)$, and

\be\label{5.9} \widetilde{Y}_t^n\leq \widetilde{Y}_t^{n+1}, 0\leq t\leq T, \mbox{a.s.},  \forall n\in {\mathbb N}.\ee

\noindent We define $\widetilde{K}_t^n:=n\int_0^t(\widetilde{{Y}}_s^n-L_s)^{-}ds,\ 0\leq t\leq T$. Then we can prove that

\bt Under our assumptions, we have\be\label{5.10}\begin{array}{lll}
&{\rm(i)}\ \widetilde{Y}_t^n\uparrow {Y}_t, 0\leq t\leq T, \mbox{a.s.};\\
&{\rm(ii)}\ E\(\sup_{0\leq t\leq T}|\widetilde{Y}_t^n-{Y}_t|^2+\int_0^T|\widetilde{Z}_s^n-{Z}_s|^2ds+\sup_{0\leq t\leq T}|\widetilde{K}_t^n-{K}_t|^2\)\rightarrow 0,
\end{array}
\ee
where $(Y, Z, K)$\ is the solution of reflected MFBSDE (\ref{5.1}).
\et
\begin{proof} We define ${f}_n(s, y, z)=E[g(s, y, z, \widetilde{Y}^n_s)],\ n\in {\mathbb N}.$ Then MFBSDE (\ref{5.8}) becomes a classical BSDE with the generator ${f}_n(s, y, z)+n(y-L_s)^{-}$\ and the terminal condition $\xi$. Here in order to be clear, we denote $C_0$\ the Lipschitz constant of $g$. Applying Lemma 2.2 to the BSDEs (\ref{5.2}) and (\ref{5.8}), for $\delta=\frac{1}{1+4C_0^2}$, there exists a constant $\gamma$\ such that
\be\begin{array}{cl}
  & \displaystyle E\[\int_0^T e^{\gamma t}(\vert{Y}^n_t-\widetilde{Y}_t^n\vert^2+\vert{Z}^n_t-\widetilde{Z}_t^n\vert^2)dt\]  \le\displaystyle  \delta E\[\int_0^T e^{\gamma t}\vert f(t, Y_t^n, Z_t^n)- {f}_n(t, Y_t^n, Z_t^n)\vert^2dt\]\\
  &\displaystyle=\delta E\[\int_0^T e^{\gamma t}\vert E[g(t, y, z, {Y}_t)]|_{y=Y_t^n, z=Z_t^n}- E[g(t, y, z, \widetilde{Y}^n_t)]|_{y=Y_t^n, z=Z_t^n}\vert^2dt\]\\
  &\displaystyle\leq 2\delta C_0^2\int_0^T e^{\gamma t}E\vert{Y}_t-\widetilde{Y}_t^n\vert^2dt\\
  &\displaystyle\leq 4\delta C_0^2\int_0^T e^{\gamma t}E\vert{Y}_t-{Y}_t^n\vert^2dt + 4\delta C_0^2\int_0^T e^{\gamma t}E\vert{Y}^n_t-\widetilde{{Y}}_t^n\vert^2dt.\\
    \end{array}\ee
Therefore, we get
\be
\displaystyle E[\int_0^T e^{\gamma t}(\vert{Y}^n_t-\widetilde{Y}_t^n\vert^2+\vert{Z}^n_t-\widetilde{Z}_t^n\vert^2)dt]  \le\displaystyle 4 C_0^2E\[\int_0^T e^{\gamma t}\vert{Y}_t-{Y}_t^n\vert^2dt\].
\ee
Furthermore, from (\ref{5.7}) and (\ref{5.9}) the proof is complete.
\end{proof}

\section{ Relation between a Reflected MFBSDE and an obstacle problem for a nonlinear
parabolic nonlocal PDE}

In this section we will show that reflected MFBSDEs studied before allow to give a probabilistic representation for the solutions of non-local PDEs with obstacles.

We consider measurable functions $b: [0,T]\times {\mathbb{R}}^n\times {\mathbb{R}}^n\rightarrow
{\mathbb{R}}^n \ $ and
         $\sigma:[0,T]\times  {\mathbb{R}}^n\times
         {\mathbb{R}}^n\rightarrow {\mathbb{R}}^{n\times d}$
which are assumed to satisfy the following conditions:
 $$
  \begin{array}{ll}
\mbox{(i)}&b(\cdot,0,0)\ \mbox{and}\ \sigma(\cdot,0, 0)\ \mbox{are continuous}
\  \mbox{and there exists}\ \mbox{some constant}\ C>0\  \mbox{such that}\\
 &\hskip0.6cm|b(t,x,\widetilde{x})|+|\sigma(t,x,\widetilde{x})|\leq C(1+|x|),\
                                  \mbox{for all}\ 0\leq t\leq T,\ x, \widetilde{x}\in {\mathbb{R}}^n;\\
\mbox{(ii)}&b\ \mbox{and}\ \sigma\ \mbox{are Lipschitz in}\ x,\ \widetilde{x},\ \mbox{i.e., there exists some constant}\ C>0\ \mbox{such that}\\
           &\hskip0.6cm|b(t,x_1,\widetilde{x}_1)-b(t, x_2, \widetilde{x}_2)|+|\sigma(t,x_1,\widetilde{x}_1)-\sigma(t, x_2, \widetilde{x}_2)|\leq C(| x_1-x_2|+| \widetilde{x}_1-\widetilde{x}_2|), \\
 & \hbox{ \ \ }\hskip7cm\mbox{for all}\ 0\leq t \leq T,\ x_1,\widetilde{x}_1, \ x_2, \widetilde{x}_2\in {\mathbb{R}}^n.\\
 \end{array}
  \eqno{\mbox{(H6.1)}}
  $$\par
  We now study the following SDE with the
  initial condition $(t,\zeta)\in[0,T]\times L^2(\Omega,{\cal{F}}_t,P;{\mathbb{R}}^n)$:
  \be\label{6.1}
  \left\{
  \begin{array}{rcl}
  dX_s^{t,\zeta}&=&E[b(s, x,X_s^{0,x_0})]|_{x=X_s^{t,\zeta}}ds+E[\sigma(s,x,X_s^{0,x_0})]|_{x=X_s^{t,\zeta}} dB_s,\ s\in[t,T],\\
  X_t^{t,\zeta}&=&\zeta.
  \end{array}
  \right.
  \ee

Under the assumption (H6.1), SDE (\ref{6.1}) has a unique strong solution.
Indeed, we first get the existence and
uniqueness of the solution $X^{0,x_0}\in
{\cal{S}}_{{\mathbb{F}}}^2(0, T; {\mathbb{R}}^n)$\ to the McKean-Vlasov
SDE (\ref{6.1}). Once knowing $X^{0,x_0}$, SDE (\ref{6.1}) becomes a classical
equation with the coefficients
$\tilde{b}(s,x)=E[b(s,x,X_s^{0,x_0})]$\
and
$\tilde{\sigma}(s,x)=E[\sigma(s,x,X_s^{0,x_0})].$\  From standard arguments we also can have, for any $p\geq
2,$\ there exists $C_{p}\in {\mathbb{R}}$\ which only depends on the
Lipschitz and the growth constants of $b$ and $\sigma$\ such that, for all
$t\in[0,T]\ \mbox{and}\ \zeta,\zeta'\in
L^p(\Omega,{\cal{F}}_t,P;{\mathbb{R}}^n),$
 \be\label{6.2-2}
 \begin{array}{rcl}
 E[\sup\limits_{t\leq s\leq T}| X_s^{t,\zeta}-X_s^{t,\zeta'}|^p|{\cal{F}}_t]
                             &\leq& C_{p}|\zeta-\zeta'|^p, \ \ \mbox{P-a.s.},\\
  E[\sup\limits_{t\leq s\leq T}| X_s^{t,\zeta}|^p|{\cal{F}}_t]
                       &\leq& C_{p}(1+|\zeta|^p),\ \ \mbox{P-a.s.}\\
  \end{array}
\ee

These standard estimates are well-known in the classical case. More details may refer to, e.g,~\cite{BLP1}.

Let now be given two real-valued mappings $f(t,x, \widetilde{x},\widetilde{y},y,z): [0,T]\times {\mathbb{R}}^n\times {\mathbb{R}}^n\times {\mathbb{R}}\times {\mathbb{R}}\times
          {\mathbb{R}}^d \rightarrow {\mathbb{R}}$ and
$\Phi(x, \widetilde{x}): {\mathbb{R}}^n\times
{\mathbb{R}}^n\rightarrow {\mathbb{R}}$ which satisfy the following conditions:
$$
\begin{array}{ll}
\mbox{(i)}&\mbox{There exists a constant}\ C>0\ \mbox{such that}\\
          &| f(t,x_1,\widetilde{x}_1,\widetilde{y}_1,y_1,z_1)-f(t,x_2,\widetilde{x}_2,\widetilde{y}_2,y_2,z_2)| +| \Phi(x_1,\widetilde{x}_1)-\Phi(x_2,\widetilde{x}_2)|\\
          &\hskip 2cm\leq C(|x_1-x_2|+|\widetilde{x}_1-\widetilde{x}_2|+ |\widetilde{y}_1-\widetilde{y}_2|+|y_1-y_2|+|z_1-z_2|),\\
&\hskip 1cm\mbox{for all}\ 0\leq t\leq T,\ x_1,\widetilde{x}_1, x_2, \widetilde{x}_2\in
{\mathbb{R}}^n,\ y_1,\widetilde{y}_1, y_2, \widetilde{y}_2 \in {\mathbb{R}}\
 \mbox{and}\ z_1, z_2\in {\mathbb{R}}^d;\\
 \mbox{(ii)}&f\ \mbox{and}\ \Phi \ \mbox{satisfy a linear growth condition, i.e., there exists some}\ C>0\\\
    & \mbox{such that},\ \mbox{for all}\ \widetilde{x},\ x\in
    {\mathbb{R}}^n,\\
    &\hskip 2cm|f(t, x,\widetilde{x},0,0,0)| + |\Phi(x,\widetilde{x})| \leq C(1+|x|+|\widetilde{x}|);\\
\mbox{(iii)}& f(t,x,\widetilde{x},\widetilde{y},y,z)\ \mbox{is continuous in}\ t\ \mbox{for all}\ (x,\widetilde{x},\widetilde{y},y,z),;\\
\mbox{(iv)}& f(t,x,\widetilde{x},  \widetilde{y},y, z)\ \mbox{is nondecreasing with respect to}\ \widetilde{y};\\
\mbox{(v)}& E[\Phi(x,X_T^{0,x_0})]\geq h(T,x),\ \mbox{for all}\ x\in
    {\mathbb{R}}^n.
\end{array}
\eqno{\mbox{(H6.2)}} $$ We consider the following reflected BSDE:
\be\label{6.2}
\begin{array}{rcl}
& & {\rm (i)}\ Y^{t,\zeta} \in S_{\mathbf{F}}^2([t,T]), \, Z^{t,\zeta} \in
L_{\mathbf{F}}^2([t,T];{\mathbb{R}}^d)\ \mbox{and}\ K^{t,\zeta}\in
A_{\mathbf{F}}^{2,c}([0,T]);\\
& & {\rm (ii)}\ Y^{t,\zeta}_s = \displaystyle E[\Phi(x,X_T^{0,x_0})]|_{x=X_T^{t,\zeta}} + \int_s^TE[f(r,x,X_r^{0,x_0},
Y_r^{0,x_0},y,z)]|_{x=X_r^{t,\zeta}, y=Y_r^{t,\zeta}, z=Z_r^{t,\zeta}}dr\\
& &\ \ \ \ \ \hskip1cm + K^{t,\zeta}_{T}-K^{t,\zeta}_{s} -\int^T_s Z^{t,\zeta}_rdW_r;\\
& & {\rm (iii)}\ \displaystyle Y^{t,\zeta}_s \geq h(s, X_s^{t,\zeta}),\ \mbox{a.s.,\ for all}\ s\in
[0,T]; \\
& & {\rm (iv)}\ \displaystyle\int_t^T(Y^{t,\zeta}_s - h(s, X_s^{t,\zeta}))dK^{t,\zeta}_{s}=0,
\end{array}
\ee
We first consider the equation (\ref{6.2}) when $(t,\zeta)=(0,
x_0)$: We know that there exists a unique solution
$(Y^{0, x_0}, Z^{0, x_0}, K^{0, x_0})\in {\cal{S}}_{{\mathbb{F}}}^2(0,
T;{\mathbb{R}})\times{L}_{{\mathbb{F}}}^{2}(0,T;{\mathbb{R}}^{d})\times A_{\mathbf{F}}^{2,c}([0,T])$\
to the Reflected Mean-Field BSDE (\ref{6.2}). Once we get $(Y^{0, x_0}, Z^{0,
x_0}, K^{0, x_0}),$\ equation (\ref{6.2}) becomes a classical reflected BSDE whose coefficients
$\tilde{f}(s,X_s^{t,\zeta},y,z)=E[f(s,x, X_s^{0,x_0},
Y_s^{0,x_0},y,z)]|_{x=X_s^{t,\zeta}}$\ satisfies the assumptions
(A1) and (A2), and $\tilde{\Phi}(X_T^{t,\zeta})=E[\Phi(x,X_T^{0,x_0})]|_{x=X_T^{t,\zeta}}\in
L^2(\Omega, {\cal{F}}_T, P)$. Thus, from Lemma 3.1 we know that
there exists a unique solution $(Y^{t,\zeta}, Z^{t,\zeta}, K^{t,\zeta})\in
{\cal{S}}_{{\mathbb{F}}}^2(0, T;
{\mathbb{R}})\times{\cal{H}}_{{\mathbb{F}}}^{2}(0,T;{\mathbb{R}}^{d})\times A_{\mathbf{F}}^{2,c}([0,T])$\
to equation (\ref{6.2}).

 Now we introduce the random field:
\be\label{6.3} u(t,x)=Y_s^{t,x}|_{s=t},\ (t, x)\in [0, T]\times{\mathbb{R}}^n,
\ee where $Y^{t,x}$ is the solution of RBSDE (\ref{6.2}) with $x \in
{\mathbb{R}}^n$\ at the place of $\zeta\in
L^2(\Omega,{\cal{F}}_t,P;{\mathbb{R}}^n).$\\

Notice that, it is obvious that $u$\ is a deterministic function, for all $t \in
[0, T], x\in {\mathbb{R}}^n $, and as we told above: once we get $(Y^{0, x_0}, Z^{0,
x_0}, K^{0, x_0}),$\ equation (\ref{6.2}) becomes a classical reflected BSDE whose coefficients
$\tilde{f}(s,X_s^{t,\zeta},y,z)=E[f(s,x,X_s^{0,x_0}, Y_s^{0,x_0},y,z)]|_{x=X_s^{t,\zeta}}$\ satisfies the assumptions
(A1) and (A2), and $\tilde{\Phi}(X_T^{t,\zeta})=E[\Phi(x,X_T^{0,x_0})]|_{x=X_T^{t,\zeta}}\in
L^2(\Omega, {\cal{F}}_T, P)$. Therefore, from Proposition 6.1 and Theorem 3.2 in~\cite{BL}, we immediately get that
 \be
\begin{array}{ll}
\mbox{(i)}&| u(t,x)-u(t,y)| \leq C|x-y|,\ \mbox{for all}\ x, y\in {\mathbb{R}}^n;\\
\mbox{(ii)}&| u(t,x)|\leq C(1+|x|),\ \mbox{for all}\ x\in {\mathbb{R}}^n;\\
\mbox{(iii)}& u\ \mbox{is continuous in}\ t.\\
\end{array}
\ee

In this section we want to consider the following non-local PDE with an obstacle \be\label{6.5}\left
\{\begin{array}{ll}
 &\!\!\!\!\! min\{u(t,x)-h(t,x),-\frac{\partial }{\partial t} u(t,x)-Au(t,x)\\
&\hskip2cm-E[f(t,x, X_t^{0, x_0},  u(t, X_t^{0, x_0}), u(t,x),
Du(t,x).E[\sigma(t, x, X_t^{0, x_0})])]\}=0,\\
 & \hskip 4.3cm (t,x)\in [0,T)\times {\mathbb{R}}^n ,  \\
 &\!\!\!\!\!  u(T,x) =E[\Phi (x, X_T^{0, x_0})], \hskip0.5cm   x \in
 {\mathbb{R}}^n,
 \end{array}\right.
\ee with
$$ Au(t,x):=\frac{1}{2}tr(E[\sigma(t,x, X_t^{0, x_0})]
 E[\sigma(t, x, X_t^{0, x_0})]^{T}D^2u(t,x))+Du(t,x).E[b(t,
 x, X_t^{0, x_0})].$$ Here the functions $b, \sigma, f\ \mbox{and}\ \Phi$\ are
supposed to satisfy (H6.1), and (H6.2), respectively, and
$X^{0, x_0}$\ is the solution of the SDE (\ref{6.1}).
We want to prove that the value function $u(t, x)$ introduced by (\ref{6.3}) is the unique viscosity solution of equation (\ref{6.5}). Now we have to do with
nonlocal PDEs with obstacles. Furthermore, unlike~\cite{BBE} here the nonlocal
term is not produced by a diffusion process with jumps. We first recall the definition
of a viscosity solution of equation (\ref{6.5}). The reader more
interested in viscosity solutions is referred to Crandall, Ishii and
Lions~\cite{CIL}.

\bde\mbox{ } A real-valued
continuous function $u\in C_p([0,T]\times {\mathbb{R}}^n )$ is called \\
  {\rm(i)} a viscosity subsolution of equation (\ref{6.5}) if, firstly, $u(T,x) \leq E[\Phi (x, X_T^{0, x_0})],\ \mbox{for all}\ x \in
  {\mathbb{R}}^n$, and if, secondly, for all functions $\varphi \in C^3_{l,b}([0,T]\times
  {\mathbb{R}}^n)$ and $(t,x) \in [0,T) \times {\mathbb{R}}^n$ such that $u-\varphi $\ attains its
  local maximum at $(t, x)$,
     $$\begin{array}{ll}
     &\min\{u(t,x)-h(t,x),-\frac{\partial }{\partial t} \varphi(t,x) - D\varphi(t,x).E[b(t,x, X_t^{0, x_0})]\\
     &-\frac{1}{2}tr(E[\sigma(t, x, X_t^{0, x_0})]
 E[\sigma(t, x, X_t^{0, x_0})]^{T}D^2\varphi(t,x))\\
 &- E[f(t, x, X_t^{0, x_0}, u(t, X_t^{0, x_0}), u(t,x), D\varphi(t,x).E[\sigma(t, x, X_t^{0, x_0})])]\}\leq 0;\end{array}
     $$
{\rm(ii)} a viscosity supersolution of equation (\ref{6.5}) if, firstly,
$u(T,x) \geq E[\Phi (x,X_T^{0, x_0})],\ \mbox{for all}\ x \in
  {\mathbb{R}}^n$, and if, secondly, for all functions $\varphi \in C^3_{l,b}([0,T]\times
  {\mathbb{R}}^n)$ and $(t,x) \in [0,T) \times {\mathbb{R}}^n$ such that $u-\varphi $\ attains its
  local minimum at $(t, x)$,
     $$\begin{array}{ll}
     &\min\{u(t,x)-h(t,x),-\frac{\partial }{\partial t} \varphi(t,x) - D\varphi(t,x).E[b(t,x, X_t^{0, x_0})]\\
     &-\frac{1}{2}tr(E[\sigma(t, x, X_t^{0, x_0})]
 E[\sigma(t, x, X_t^{0, x_0})]^{T}D^2\varphi(t,x))\\
 &- E[f(t, x, X_t^{0, x_0}, u(t, X_t^{0, x_0}), u(t,x), D\varphi(t,x).E[\sigma(t, x, X_t^{0, x_0})])]\} \geq 0;\end{array}
     $$
 {\rm(iii)} a viscosity solution of equation (\ref{6.5}) if it is both a viscosity sub- and a supersolution of equation
     (6.1).\ede
\br {\rm(i)} We recall that $C_p([0,T]\times R^n)=\{u\in C([0,T]\times R^n):
\mbox{There exists some constant}\ p>0\ \mbox{such that}$ $
\sup_{(t,x)\in[0, T]\times R^n}\frac{|u(t,x)|}{1+|x|^p}<
+\infty\}.$\\
{\rm(ii)} The space
$C^3_{l,b}([0,T]\times {\mathbb{R}}^n)$ denotes the set of the
real-valued functions that are continuously differentiable up to the
third order and whose derivatives of order from 1 to 3 are bounded.
Therefore, that function in $C^3_{l,b}([0,T]\times
{\mathbb{R}}^n)$ is of at most linear growth.\er

We now can give the main statement of this section.
 \bt Under the assumptions (H6.1) and (H6.2) the function $u(t,x)$\ defined by (\ref{6.3}) is the unique viscosity solution of
equation (\ref{6.5}). \et

For each $(t,x)\in [0,T]\times{\mathbb{R}^n}$, and $n\in\mathbf{N}$,
let $\{(^nY^{t,x}_s,{}^nZ^{t,x}_s), t\leq s\leq T\}$ denote
the solution of the MFBSDE
\be\label{6.6}\aligned
^nY^{t,x}_s &=E[\Phi(x, X_T^{0,x_0})]|_{x=X_T^{t,x}}+\int_s^T
E[f(r,x,X_r^{0,x_0},
^n\!Y_r^{0,x_0},y,z)]|_{x=X_r^{t,x},y=^n\!Y^{t,x}_r,z=^n\!Z^{t,x}_r }dr\\
&\quad +n\int_s^T(^nY^{t,x}_r-h(r,X^{t,x}_r))^-dr-\int_s^T\
^nZ^{t,x}_rdW_r,\quad t\leq s\leq T.
\endaligned\ee
We define \be\label{6.7} u_n(t,x) := ^n\!Y^{t,x}_t,\   0\leq t\leq T,\ x\in{\mathbb{R}^n}.\ee

It is known from~\cite{BLP} that $u_n(t,x)$ defined by (\ref{6.7}) is in $C([0,T]\times
{\mathbb{R}^n})$, has linear growth in $x$, and is the unique
continuous viscosity solution of the following equation:
\be\label{6.8}\left
\{\begin{array}{ll}
 &\!\!\!\!\! -\frac{\partial }{\partial t} u_n(t,x) -Au_n(t,x)
- \{E[f(t,x, X_t^{0, x_0},  u_n(t, X_t^{0, x_0}), u_n(t,x),
Du_n(t,x).E[\sigma(t,x, X_t^{0, x_0})])]\\
 &+n(u_n(t,x)-h(t,x))^-\}=0, \hskip 4.3cm (t,x)\in [0,T)\times {\mathbb{R}}^n ,  \\
 &\!\!\!\!\!  u_n(T,x) =E[\Phi (x,X_T^{0, x_0})], \hskip0.5cm   x \in
 {\mathbb{R}}^n,
 \end{array}\right.
\ee with
$$ Au(t,x):=\frac{1}{2}tr(E[\sigma(t,x, X_t^{0, x_0})]
 E[\sigma(t,x, X_t^{0, x_0})]^{T}D^2u(t,x))+Du(t,x).E[b(t,x,
 X_t^{0, x_0})].$$
We have the uniqueness of viscosity solution $u_n$\ only in the
space $C_p([0,T]\times R^n)$ (in~\cite{BLP} the authors gave an example to explain why the uniqueness is only in $C_p([0,T]\times R^n)$). More details refer to~\cite{BLP}.

\bl
\be\label{6.8-1} u_n(t,x)\uparrow u(t,x),\ \mbox{for all}\ (t, x)\in [0, T]\times
{\mathbb{R}}^n.\ee
\el
\noindent{\bf Proof}. When $(t,x)=(0,x_0)$, the equation (\ref{6.6}) is the penalized reflected mean-field BSDE, from Section 5, we know $^n\!Y^{0,x_0}_r\uparrow Y_r^{0,x_0}, \ 0\leq r\leq T$, P-a.s., in particular, $u_n(0,x_0)\uparrow u(0,x_0)$.

When $(t,x)\neq(0,x_0)$, recall that SDE (\ref{6.1}) becomes the classical
equation with the coefficients $\tilde{b}(r,x)=E[b(r,x,X_r^{0,x_0})]$\ and $\tilde{\sigma}(r,x)=E[\sigma(r,x,X_r^{0,x_0})],$\ the equation (\ref{6.6}) becomes the classical panalized BSDE with the coefficient $\widetilde{f}_n(r,x, y,z):=E[f(r,x,X_r^{0,x_0},
^n\!Y_r^{0,x_0},y,z)]+n(y-h(r,x))^-$, and terminal condition $\widetilde{\Phi}(x):=E[\Phi (x,X_T^{0, x_0})]$. Notice that now still $\widetilde{f}_n(r,x, y,z)\leq \widetilde{f}_{n+1}(r,x, y,z)$, following the proof on Pages 719-723 in~\cite{EKPPQ}, we can get $^n\!Y^{t,x}_s\uparrow Y_s^{t,x}, \ t\leq s\leq T$, P-a.s., therefore, $u_n(t,x)\uparrow u(t,x),\ \mbox{for all}\ (t, x)\in [0, T]\times
{\mathbb{R}}^n.$\endpf

On the other hand, notice that because $u_n$ and $u$\ are continuous, from Dini's theorem it follows that the above convergence is uniform on compacts.
We can also prove that, $u_n$\ has linear growth in $x$\ and is Lipschitz in $x$, uniformly with respect to $n\in \mathbf{N}$.
\bp There exists a constant $C$\ independent of $n$, such that, for every $n\in \mathbf{N}$,
\be\label{6.10}\begin{array}{lll}
& {\rm (i)}\ |u_n(t,x)|\leq C(1+|x|), \ \mbox{for all}\ x \in {\mathbb{R}}^n,\ t\in [0, T];\\
& {\rm (ii)}\ | u_n(t,x)-u_n(t,y)| \leq C|x-y|,\ \mbox{for all}\ x, y\in {\mathbb{R}}^n,\ t\in [0, T].\end{array}\ee
\ep
The proof is given in the appendix for convenience.

\medskip

\noindent{\bf Proof of Theorem 6.1}. \textbf{Step 1}: We first prove that $u(t,x)$\ is a viscosity supersolution of (\ref{6.5}).
Indeed, let $(t,x)\in [0,T)\times {\mathbb{R}}^n$ and
let $\varphi\in C^3_{l, b}([0,T]\times {\mathbb{R}}^n)$\ be such
that $u-\varphi>u(t,x)-\varphi(t,x)=0$\
everywhere on $([0,T]\times {\mathbb{R}}^n) -\{(t,x)\}.$\ Then,
because $u$\ is continuous and $u_n(t,x)\uparrow
u(t,x)$, $0\leq t\leq T $, $x\in{\mathbb{R}}^n$, there
exists some sequence $(t_n,x_n),\ n\geq 1,$\ at least
along a subsequence, such that,

\smallskip

i) $(t_n,x_n)\rightarrow (t,x)$, as $n\rightarrow +\infty$;

\smallskip
ii) $u_n-\varphi\geq u_n(t_n,x_n)- \varphi(t_n,x_n)$ in a
neighborhood of $(t_n,x_n)$, for all $n\geq 1$;

\smallskip
iii) $u_n(t_n,x_n)\rightarrow u(t,x)$, as $n\rightarrow
+\infty$.

\medskip

\noindent Consequently, because $u_n$ is a viscosity solution and
hence a supersolution of equation (\ref{6.8}), we have, for all $n\geq 1$,
\be\begin{array}{ll}
 &\!\!\!\!\! \frac{\partial }{\partial t}\varphi(t_n,x_n)+\frac{1}{2}tr(E[\sigma(t_n, x_n, X_{t_n}^{0, x_0})]
 E[\sigma(t_n, x_n, X_{t_n}^{0, x_0})]^{T}D^2\varphi(t_n,x_n))+D\varphi(t_n,x_n).E[b(t_n,
 x_n, X_{t_n}^{0, x_0})]\\
&+\{E[f(t_n, x_n, X_{t_n}^{0, x_0}, u_n(t_n, X_{t_n}^{0, x_0}), u_n(t_n,x_n),
D\varphi(t_n,x_n).E[\sigma(t_n, x_n, X_{t_n}^{0, x_0})])]\\
 &+n(u_n(t_n,x_n)-h(t_n,x_n))^-\}\leq 0,
 \end{array}
\ee
\noindent Therefore,
\be
\begin{array}{ll}
 &\!\!\!\!\! \frac{\partial }{\partial t}\varphi(t_n,x_n)+\frac{1}{2}tr(E[\sigma(t_n, x_n, X_{t_n}^{0, x_0})]
 E[\sigma(t_n, x_n, X_{t_n}^{0, x_0})]^{T}D^2\varphi(t_n,x_n))+D\varphi(t_n,x_n).E[b(t_n,
 x_n, X_{t_n}^{0, x_0})]\\
&+E[f(t_n, x_n, X_{t_n}^{0, x_0}, u_n(t_n, X_{t_n}^{0, x_0}), u_n(t_n,x_n),
D\varphi(t_n,x_n).E[\sigma(t_n, x_n, X_{t_n}^{0, x_0})])]\leq 0,
 \end{array}
\ee
Taking the limit, from (H6.1), (H6.2) and (\ref{6.10}) it follows from Lebesgue dominated convergence theorem we get that
\be
\begin{array}{ll}
 &\!\!\!\!\! \frac{\partial }{\partial t}\varphi(t,x)+\frac{1}{2}tr(E[\sigma(t,x, X_{t}^{0, x_0})]
 E[\sigma(t,x, X_{t}^{0, x_0})]^{T}D^2\varphi(t,x))+D\varphi(t,x).E[b(t,x,
 X_{t}^{0, x_0})]\\
&+E[f(t, x,X_{t}^{0, x_0}, u(t, X_{t}^{0, x_0}), u(t,x),
D\varphi(t,x).E[\sigma(t,x, X_{t}^{0, x_0})])]\leq 0,
 \end{array}
\ee
Because our $u(t,x)\geq h(t,x)$, we prove that $u(t,x)$ is the viscosity supersolution.

\textbf{Step 2}: The function $W(t,x)$\ is a viscosity subsolution of  equations (\ref{6.5}).

Indeed, let $(t,x)\in [0,T)\times {\mathbb{R}}^n$\ be a point at which $u(t,x)>h(t,x)$, and
let $\varphi\in C^3_{l, b}([0,T]\times {\mathbb{R}}^n)$\ be such
that $u-\varphi<u(t,x)-\varphi(t,x)=0$\
everywhere on $([0,T]\times {\mathbb{R}}^n) -\{(t,x)\}.$\ Then,
because $u$\ is continuous and $u_n(t,x)\uparrow
u(t,x)$, $0\leq t\leq T $, $x\in{\mathbb{R}}^n$, there
exists some sequence $(t_n,x_n),\ n\geq 1,$\ at least
along a subsequence, such that,

\smallskip

i)$(t_n,x_n)\rightarrow (t,x)$, as $n\rightarrow +\infty$;

\smallskip
ii) $u_n-\varphi\leq u_n(t_n,x_n)- \varphi(t_n,x_n)$ in a
neighborhood of $(t_n,x_n)$, for all $n\geq 1$;

\smallskip
iii) $u_n(t_n,x_n)\rightarrow u(t,x)$, as $n\rightarrow
+\infty$.

\medskip

\noindent Consequently, because $u_n$ is a viscosity solution and
hence a subsolution of equation (\ref{6.8}), we have, for all $n\geq 1$,
\be\begin{array}{ll}
 &\!\!\!\!\! \frac{\partial }{\partial t}\varphi(t_n,x_n)+\frac{1}{2}tr(E[\sigma(t_n, x_n, X_{t_n}^{0, x_0})]
 E[\sigma(t_n, x_n, X_{t_n}^{0, x_0})]^{T}D^2\varphi(t_n,x_n))+D\varphi(t_n,x_n).E[b(t_n,
 x_n, X_{t_n}^{0, x_0})]\\
&+\{E[f(t_n, x_n, X_{t_n}^{0, x_0}, u_n(t_n, X_{t_n}^{0, x_0}), u_n(t_n,x_n),
D\varphi(t_n,x_n).E[\sigma(t_n, x_n, X_{t_n}^{0, x_0})])]\\
 &+n(u_n(t_n,x_n)-h(t_n,x_n))^-\}\geq 0,
 \end{array}
\ee

From the assumption that $u(t,x)>h(t,x)$ and the uniform convergence of $u_n$, we can get that for $n$\ large enough $u_n(t_n,x_n)>h(t_n,x_n)$, hence, taking the limit as $n\rightarrow \infty$ in the above inequality from (H6.1), (H6.2) and (\ref{6.10}) it follows from Lebesgue dominated convergence theorem we get:

\be\begin{array}{ll}
 &\!\!\!\!\! \frac{\partial }{\partial t}\varphi(t,x)+\frac{1}{2}tr(E[\sigma(t,x, X_{t}^{0, x_0})]
 E[\sigma(t, x, X_{t}^{0, x_0})]^{T}D^2\varphi(t,x))+D\varphi(t,x).E[b(t,
 x, X_{t}^{0, x_0})]\\
&+E[f(t, x, X_{t}^{0, x_0}, u(t, X_{t}^{0, x_0}), u(t,x),
D\varphi(t,x).E[\sigma(t, x, X_{t}^{0, x_0})])]\\
 &\geq 0,
 \end{array}
\ee we prove that $u(t,x)$ is the viscosity subsolution.
\endpf

\br For the uniqueness of the viscosity solution, from~\cite{BLP} we know that we have the uniqueness of viscosity solution $u$\ in the
space $C_p([0,T]\times R^n)$. Similar to the proof of Theorem 5.1 in~\cite{BL} and Theorem 7.1 in~\cite{BLP} it is not hard to prove the uniqueness.
\er

\section{Appendix}

\noindent{\bf Proof of Proposition 6.1}. (i) From (\ref{6.8-1}) we know $u_1(t,x)\leq u_n(t,x)\leq u(t,x)$, for all $n\in {\mathbb{N}}$, \ $x \in {\mathbb{R}}^n,\ t\in [0, T].$\ Furthermore, since $u_1$\ and $u$ are linear growth in $x$, the proof of (i) is complete.

(ii) To simplify the notations, we define $\tilde{b}(s,x)=E[b(s,x,X_s^{0,x_0})]$\
and $\tilde{\sigma}(s,x)=E[\sigma(s,x,X_s^{0,x_0})],$ $\widetilde{\Phi}(x)=E[\Phi(x,X_T^{0,x_0})], \widetilde{f}(s, x, y, z)=E[f(s,x,X_s^{0,x_0},
^n\!Y_s^{0,x_0},y,z)]. $
Then, from the equation (\ref{6.6}) we know $(^nY^{t,x}, ^n\!Z^{t,x})$ is the solution of the following equation:
\be\label{100}\aligned
^nY^{t,x}_s &=\widetilde{\Phi}(X_T^{t,x})+\int_s^T
\widetilde{f}(r,X_r^{t,x},^n\!Y^{t,x}_r, ^n\!Z^{t,x}_r)dr\\
&\quad +n\int_s^T(^nY^{t,x}_r-h(r,X^{t,x}_r))^-dr-\int_s^T\
^nZ^{t,x}_rdW_r,\quad t\leq s\leq T,
\endaligned\ee
and  $(^nY^{t,y}, ^n\!Z^{t,y})$ is the solution of the following equation:
\be\label{101}\aligned
^nY^{t,y}_s &=\widetilde{\Phi}(X_T^{t,y})+\int_s^T
\widetilde{f}(r,X_r^{t,y},^n\!Y^{t,y}_r, ^n\!Z^{t,y}_r)dr\\
&\quad +n\int_s^T(^nY^{t,y}_r-h(r,X^{t,y}_r))^-dr-\int_s^T\
^nZ^{t,y}_rdW_r,\quad t\leq s\leq T.
\endaligned\ee
For an arbitrarily given $\varepsilon>0$, we consider the function $\psi_\varepsilon(x)=(|x|^2+\varepsilon)^{\frac{1}{2}},\ x\in
 {\mathbb{R}}^n.$\ Obviously, $|x|\leq \psi_\varepsilon(x)\leq |x|+\varepsilon^{\frac{1}{2}},\ x\in
 {\mathbb{R}}^n.$\ Moreover, for all $x\in {\mathbb{R}}^n,$
 $$
D\psi_\varepsilon(x)=\frac{x}{(|x|^2+\varepsilon)^{\frac{1}{2}}},\ \
\ \ \
D^2\psi_\varepsilon(x)=\frac{I}{(|x|^2+\varepsilon)^{\frac{1}{2}}}-\frac{x\otimes
x}{(|x|^2+\varepsilon)^{\frac{3}{2}}}.
 $$
Therefore, we have \be\label{102} |D\psi_\varepsilon(x)|\leq 1,\ \ \
|D^2\psi_\varepsilon(x)||x|\leq
\frac{C}{(|x|^2+\varepsilon)^{\frac{1}{2}}}|x|\leq C,\ \ x\in
{\mathbb{R}}^n,\ee where the constant $C$\ is independent of
$\varepsilon$. On the other hand, in order to be clear we denote $\mu$\ is the
Lipschitz constant of $h, \ \Phi,$ and $f$. We consider the
following two BSDEs: \be\label{103}
\begin{array}{lll}
 &\tilde{Y}_s = \widetilde{\Phi}(X_T^{t,x}) +
\mu\psi_\varepsilon(X_T^{t,x}-X_T^{t,y})+
\int_s^T(\widetilde{f}(r,X_r^{t,x},\tilde{Y}_r,\tilde{Z}_r)+\mu|X_r^{t,x}-X_r^{t,y}|)dr\\
&\ \hskip1cm + \int_s^Tn\(\tilde{Y}_r-\(h(r,
X_r^{t,x})+\mu\psi_\varepsilon(X_r^{t,x}-X_r^{t,y})\)\)^{-}dr- \int^T_s\tilde{Z}_rdW_r,\ \ \ \  s\in [t,T];\ \\
\end{array}
\ee

\noindent and
\be\label{104}
\begin{array}{lll}
 & \bar{Y}_s = \widetilde{\Phi}(X_T^{t,x}) -
\mu\psi_\varepsilon(X_T^{t,x}-X_T^{t,y})+
\int_s^T(\widetilde{f}(r,X_r^{t,x},\bar{Y}_r,\bar{Z}_r)-\mu|X_r^{t,x}-X_r^{t,y}|)dr\\
&\ \hskip1cm + \int_s^Tn\(\bar{Y}_r-\(h(r,
X_r^{t,x})-\mu\psi_\varepsilon(X_r^{t,x}-X_r^{t,y})\)\)^{-}dr- \int^T_s\bar{Z}_rdW_r,\ \ \ \  s\in [t,T].
\end{array}\ee

\noindent Obviously, the coefficients satisfy the assumptions (A1) and (A2), therefore from Lemma 2.1 they have unique solutions $(\tilde{Y}, \tilde{Z})$\ and $(\bar{Y}, \bar{Z})$, respectively. Notice that the solutions of (\ref{103}) and (\ref{104}) depend on $n$, for simplifying notations and causing no confusion we still denote the solutions by $(\tilde{Y}, \tilde{Z})$\ and $(\bar{Y}, \bar{Z})$, respectively.
Furthermore, from the comparison theorem for BSDEs (Lemma 2.3)
\be\label{104-1}\bar{Y}_s\leq ^n\!Y^{t,x}_s\leq \tilde{Y}_s,\ \ \ \bar{Y}_s\leq
^n\!Y^{t,y}_s\leq \tilde{Y}_s,\ \ \mbox{P-a.s., for all}\ s\in [t,
T].\ee Now we introduce two other BSDEs: \be\label{105}
\begin{array}{lll}
 &\tilde{Y}'_s = \widetilde{\Phi}(X_T^{t,x}) +\\
&\ \ \ \
\int_s^T[\widetilde{f}(r,X_r^{t,x},\tilde{Y}'_r+\mu\psi_\varepsilon(X_r^{t,x}-X_r^{t,y}),\tilde{Z}'_r
+\mu D\psi_\varepsilon(X_r^{t,x}-X_r^{t,y})(\sigma(r,X_r^{t,x})-\sigma(r,X_r^{t,y})))\\
&\ \hskip1cm +\mu|X_r^{t,x}-X_r^{t,y}|+n(\tilde{Y}'_r-h(r,X_r^{t,x}))^{-}+\mu D\psi_\varepsilon(X_r^{t,x}-X_r^{t,y})(\widetilde{b}(r,X_r^{t,x})-\widetilde{b}(r,X_r^{t,y}))\\
&\ \hskip1cm+\frac{1}{2}\mu(D^2\psi_\varepsilon(X_r^{t,x}-X_r^{t,y})
(\widetilde{\sigma}(r,X_r^{t,x})-\widetilde{\sigma}(r,X_r^{t,y})),\widetilde{\sigma}(r,X_r^{t,x})-\widetilde{\sigma}(r,X_r^{t,y}))]dr\\
&\ \ \ \  - \int^T_s\tilde{Z}'_rdW_r,\ \ \ \  s\in [t,T];\ \\
\end{array} \ee
\noindent and \be\label{106}
\begin{array}{lll}
 &\bar{Y}'_s = \widetilde{\Phi}(X_T^{t,x})\\
&\ \ \ \ \int_s^T[\widetilde{f}(r,X_r^{t,x},\bar{Y}'_r-\mu\psi_\varepsilon(X_r^{t,x}-X_r^{t,y}),\bar{Z}'_r
-\mu D\psi_\varepsilon(X_r^{t,x}-X_r^{t,y})(\sigma(r,X_r^{t,x})-\sigma(r,X_r^{t,y})))\\
&\ \hskip1cm -\mu|X_r^{t,x}-X_r^{t,y}|+n(\bar{Y}'_r-h(r,X_r^{t,x}))^{-}-\mu D\psi_\varepsilon(X_r^{t,x}-X_r^{t,y})(\widetilde{b}(r,X_r^{t,x})-\widetilde{b}(r,X_r^{t,y}))\\
&\
\hskip1cm-\frac{1}{2}\mu(D^2\psi_\varepsilon(X_r^{t,x}-X_r^{t,y})
(\widetilde{\sigma}(r,X_r^{t,x})-\widetilde{\sigma}(r,X_r^{t,y})),\widetilde{\sigma}(r,X_r^{t,x})-\widetilde{\sigma}(r,X_r^{t,y}))]dr\\
&\ \ \ \ - \int^T_s\bar{Z}'_rdW_r,\ \ \ \  s\in [t,T].\ \\
\end{array} \ee
It's obvious that the coefficients of the BSDEs (\ref{105}) and (\ref{106}) satisfy the assumptions (A1) and (A2), hence, from Lemma 2.1 (\ref{105}) and (\ref{106}) have unique solutions $(\tilde{Y}', \tilde{Z}')$\ and $(\bar{Y}', \bar{Z}')$, respectively. On
the other hand, from the uniqueness of the solution of BSDE we know
that \be\label{107}\begin{array}{lll}
&\widetilde{Y}'_s=\widetilde{Y}_s-\mu\psi_\varepsilon(X_s^{t,x}-X_s^{t,y}),\
\mbox{for all}\ s\in [t, T],\ \mbox{P-a.s.,}\\
&\tilde{Z}'_s=\tilde{Z}_s-\mu
D\psi_\varepsilon(X_s^{t,x}-X_s^{t,y})(\widetilde{\sigma}(s,X_s^{t,x}
)-\widetilde{\sigma}(s,X_s^{t,y})),\  \mbox{dsdP-a.e. on}\ [t,
T]\times\Omega;\end{array}\ee \noindent and \be\label{107-1}\begin{array}{lll}
&\bar{Y}'_s=\bar{Y}_s+\mu\psi_\varepsilon(X_s^{t,x}-X_s^{t,y}),\
\mbox{for all}\ s\in [t, T],\ \mbox{P-a.s.,}\\
&\bar{Z}'_s=\bar{Z}_s+\mu
D\psi_\varepsilon(X_s^{t,x}-X_s^{t,y})(\widetilde{\sigma}(s,X_s^{t,x}
)-\widetilde{\sigma}(s,X_s^{t,y})),\  \mbox{dsdP-a.e. on}\ [t,
T]\times\Omega.\end{array}\ee

For any $\beta>0$ applying It\^{o}'s formula to $e^{\beta s}|\widetilde{Y}'_s-\bar{Y}'_s|^2$ we get
\be\label{108}\begin{array}{lll}
&e^{\beta r}|\widetilde{Y}'_r-\bar{Y}'_r|^2+E[\int_r^T\beta e^{\beta s}|\widetilde{Y}'_s-\bar{Y}'_s|^2ds|{\cal F}_r]+E[\int_r^Te^{\beta s}|\widetilde{Z}'_s-\bar{Z}'_s|^2ds|{\cal F}_r]=\\
&E[\int_r^T2e^{\beta s}(\widetilde{Y}'_s-\bar{Y}'_s)(\widetilde{f}(s,X_s^{t,x},\tilde{Y}'_s+\mu\psi_\varepsilon(X_s^{t,x}-X_s^{t,y}),\tilde{Z}'_s
+\mu D\psi_\varepsilon(X_s^{t,x}-X_s^{t,y})(\widetilde{\sigma}(s,X_s^{t,x})-\widetilde{\sigma}(s,X_s^{t,y})))\\
&-\widetilde{f}(s,X_s^{t,x},\bar{Y}'_s-\mu\psi_\varepsilon(X_s^{t,x}-X_s^{t,y}),\bar{Z}'_s
-\mu D\psi_\varepsilon(X_s^{t,x}-X_s^{t,y})(\widetilde{\sigma}(s,X_s^{t,x})-\widetilde{\sigma}(s,X_s^{t,y}))))ds|{\cal F}_r]\\
&+E[\int_r^T2 e^{\beta s}n(\widetilde{Y}'_s-\bar{Y}'_s)((\tilde{Y}'_s-h(s,X_s^{t,x}))^{-}-(\bar{Y}'_s-h(s,X_s^{t,x}))^{-})ds|{\cal F}_r]\\
&+E[\int_r^T2 e^{\beta s}(\widetilde{Y}'_s-\bar{Y}'_s)\Delta g(s)ds|{\cal F}_r],\\
\end{array}\ee
\noindent where
\be\begin{array}{lll} &\Delta g(s)= 2\mu|X_s^{t,x}-X_s^{t,y}|+2\mu D\psi_\varepsilon(X_s^{t,x}-X_s^{t,y})(\widetilde{b}(s,X_s^{t,x})-\widetilde{b}(s,X_s^{t,y}))\\
&\ \hskip1.5cm+\mu(D^2\psi_\varepsilon(X_s^{t,x}-X_s^{t,y})
(\widetilde{\sigma}(s,X_s^{t,x})-\widetilde{\sigma}(s,X_s^{t,y})),\widetilde{\sigma}(s,X_s^{t,x})-\widetilde{\sigma}(s,X_s^{t,y})).\\
\end{array}\ee
Notice that $(\widetilde{Y}'_s-\bar{Y}'_s)((\tilde{Y}'_s-h(s,X_s^{t,x}))^{-}-(\bar{Y}'_s-h(s,X_s^{t,x}))^{-})\leq 0.$

From (\ref{102}) and the Lipschitz continuity of $\widetilde{b}$\ and $\widetilde{\sigma}$\
we get $|\Delta g(s)|\leq C|X_s^{t,x}-X_s^{t,y}|,\ \mbox{P-a.s.},$\ where the constant $C$ \ is independent
of $\varepsilon$\ and $n$. Therefore, from (\ref{108}) we have
\be\label{109}\begin{array}{lll}
&e^{\beta r}|\widetilde{Y}'_r-\bar{Y}'_r|^2+E[\int_r^T\beta e^{\beta s}|\widetilde{Y}'_s-\bar{Y}'_s|^2ds|{\cal F}_r]+E[\int_r^Te^{\beta s}|\widetilde{Z}'_s-\bar{Z}'_s|^2ds|{\cal F}_r]\\
\leq &C E[\int_r^Te^{\beta s}|\widetilde{Y}'_s-\bar{Y}'_s|^2ds|{\cal F}_r]+\frac{1}{2}E[\int_r^T2 e^{\beta s}|\widetilde{Z}'_s-\bar{Z}'_s|^2ds|{\cal F}_r]\\
&+CE[\int_r^T e^{\beta s}|X_s^{t,x}-X_s^{t,y}|^2ds|{\cal F}_r]+C\varepsilon,\ \mbox{P-a.s.}, r\in [t, T],\\
\end{array}\ee
where the constant $C$ \ is independent
of $\varepsilon$\ and $n$. Then, take $\beta=C+1$ with the help of (\ref{6.2-2}) we get
\be\label{110}\begin{array}{lll}
|\widetilde{Y}'_t-\bar{Y}'_t|^2&\leq& CE[\int_t^T|X_s^{t,x}-X_s^{t,y}|^2ds|{\cal F}_t]+C\varepsilon\\
&\leq& CE[\sup_{t\leq s\leq T}|X_s^{t,x}-X_s^{t,y}|^2|{\cal F}_t]+C\varepsilon\\
&\leq& C|x-y|^2+C\varepsilon,\  \mbox{P-a.s.}
\end{array}\ee
Furthermore, from (\ref{6.7}), (\ref{104-1}), (\ref{107}), (\ref{107-1}) and (\ref{110}) we have
$$\begin{array}{rcl}
& &|u_n(t,x)-u_n(t,y)|^2 =|^n\!Y^{t,x}_t-^n\!Y^{t,y}_t|^2\\
& &\leq |\tilde{Y}_t-\bar{Y}_t|^2 \leq 2|\tilde{Y}'_t-\bar{Y}'_t|^2+16\mu^2(|X^{t,x}_t-X^{t,y}_t|^2+\varepsilon)\\
& &\leq  C|x-y|^2+C\varepsilon.\end{array}$$

\noindent Then, let $\varepsilon$\ tend to 0 the proof is complete.\endpf

\end{document}